\documentclass[final,onefignum,onetabnum]{siamart190516}

\usepackage{lipsum}
\usepackage{amsfonts}
\usepackage{graphicx}
\usepackage{epstopdf}
\usepackage{algorithmic}
\usepackage{amsopn}

\usepackage{amsmath}
\usepackage{amssymb}
\usepackage{fancyhdr}
\usepackage{color}
\usepackage[caption=false]{subfig}
\usepackage{multirow}
\usepackage{rotating}
\usepackage{tabularx}
 \allowdisplaybreaks[4]

\newcommand{\tr}{^{\sf T}}
\newsiamremark{remark}{Remark}
\newsiamremark{hypothesis}{Hypothesis}
\crefname{hypothesis}{Hypothesis}{Hypotheses}
\newsiamthm{claim}{Claim}
\title{On the asymptotic convergence and acceleration of gradient methods\thanks{August 19, 2019,
This research was supported by the National Natural Science Foundation of
China (11701137, 11631013, 11671116), by the National 973 Program of China
(2015CB856002), by the China Scholarship Council (No. 201806705007), and by the
USA National Science Foundation (1522654, 1819161).}}

\author{Yakui Huang\thanks{Institute of Mathematics, Hebei University of Technology, Tianjin 300401, China
  (\email{huangyakui2006@gmail.com}).}
\and Yu-Hong Dai\thanks{LSEC, Academy of Mathematics and Systems Science, Chinese Academy of Sciences, Beijing 100190, China
  (\email{dyh@lsec.cc.ac.cn}, \url{http://lsec.cc.ac.cn/\string~dyh/}).}
\and Xin-Wei Liu\thanks{Institute of Mathematics, Hebei University of Technology, Tianjin 300401, China
  (\email{mathlxw@hebut.edu.cn}).}
\and Hongchao Zhang\thanks{Department of Mathematics, Louisiana State University, Baton Rouge, LA 70803-4918, USA
(\email{hozhang@math.lsu.edu}, \url{https://www.math.lsu.edu/\string~hozhang/}).}
}

\begin{document}

\maketitle

% REQUIRED
\begin{abstract}
  We consider the asymptotic behavior of a family of gradient methods, which include the steepest descent and minimal gradient methods as special instances. It is proved that each method in the family will asymptotically zigzag between two directions. Asymptotic convergence results of the objective value, gradient norm, and stepsize are presented as well. To accelerate the family of gradient methods,
we further exploit spectral properties of stepsizes to break the zigzagging pattern. In particular, a new stepsize is derived by imposing finite termination on minimizing two-dimensional strictly convex
quadratic function. It is shown that, for the general quadratic function, the proposed stepsize asymptotically converges to the reciprocal of the largest eigenvalue of the Hessian.
Furthermore, based on this spectral property, we propose a periodic gradient method by incorporating the Barzilai-Borwein method. Numerical comparisons with some recent successful gradient methods show that our new method is very promising.
\end{abstract}

% REQUIRED

% REQUIRED

\pagestyle{myheadings}
\thispagestyle{plain}
\markboth{Y. Huang, Y-H, Dai, X-W, Liu and H. ZHANG}
{ASYMPTOTIC CONVERGENCE AND ACCELERATION OF GRADIENT METHODS}

\section{Introduction}
\label{intro}

The gradient method is well-known for solving the following unconstrained optimization
\begin{equation}\label{eqpro}
  \min_{x\in\mathbb{R}^n}~f(x),
\end{equation}
where $f: \mathbb{R}^n \to \mathbb{R}$ is continuously differentiable, especially when the dimension $n$ is large.
In particular, at $k$-th iteration gradient methods update the iterates by
\begin{equation}\label{eqitr}
  x_{k+1}=x_k-\alpha_kg_k,
\end{equation}
where $g_k=\nabla f(x_k)$ and $\alpha_k>0$ is the stepsize determined by the method.

One simplest nontrivial nonlinear instance of \eqref{eqpro} is the quadratic optimization
\begin{equation}\label{qudpro}
  \min_{x\in\mathbb{R}^n}~f(x)=\frac{1}{2}x \tr Ax-b \tr x,
\end{equation}
where $b\in\mathbb{R}^n$ and $A\in\mathbb{R}^{n\times n}$ is symmetric and positive definite.
Solving \eqref{qudpro} efficiently is usually a pre-requisite for a method to be generalized to solve more general optimization. In addition, by Taylor's expansion, a general smooth function
can be approximated by a quadratic function near the minimizer. So, the local convergence behaviors of gradient methods are often reflected by solving \eqref{qudpro}.
Hence, in this paper, we focus on studying the convergence behaviors and propose efficient gradient methods for solving \eqref{qudpro} efficiently.

In \cite{cauchy1847methode}, Cauchy proposed the steepest descent (SD) method that solves \eqref{qudpro}  by using the exact stepsize
\begin{equation}\label{eqsd}
  \alpha_k^{SD}=\arg\min_{\alpha}~f(x_k-\alpha g_k)=\frac{g_{k} \tr g_{k}}{g_{k} \tr Ag_{k}}.
\end{equation}
Although $\alpha_k^{SD}$ minimizes $f$ along the steepest descent direction, the SD method often performs poorly in practice
and has linear converge rate  \cite{akaike1959successive,forsythe1968asymptotic} as
\begin{equation}\label{sdrate}
  \frac{f(x_{k+1})-f^*}{f(x_{k})-f^*}\leq\left(\frac{\kappa-1}{\kappa+1}\right)^2,
\end{equation}
where $f^*$ is the optimal function value of \eqref{qudpro} and $\kappa=\lambda_n/\lambda_1$ is the condition number of $A$ with $\lambda_1$ and $\lambda_n$ being the smallest and largest eigenvalues of $A$, respectively.
%$R_{SD}$ is the worst asymptotic rate given by
%\begin{equation}\label{ratesd}
%  R_{SD}=\left(\frac{\kappa-1}{\kappa+1}\right)^2,
%\end{equation}
%with $\kappa=\frac{\lambda_n}{\lambda_1}$ being the condition number of $A$. Here, $\lambda_1$ and $\lambda_n$ are the smallest and largest eigenvalues of $A$, respectively.
Thus, if $\kappa$ is large, the SD method may converge very slowly. In addition, Akaike \cite{akaike1959successive} proved that the gradients will asymptotically alternate between two directions in the subspace spanned by the two eigenvectors corresponding to $\lambda_1$ and $\lambda_n$. So, the SD method often has zigzag phenomenon near the solution.
In \cite{forsythe1968asymptotic}, Forsythe generalized Akaike's results to the so-called optimum $s$-gradient method and Pronzato et al. \cite{pronzato2006asymptotic} further generalized the results
to the so-called $P$-gradient methods in the Hilbert space.
Recently, by employing Akaike's results, Nocedal et al. \cite{nocedal2002behavior} presented some insights for asymptotic behaviors of the SD method on function values, stepsizes and gradient norms.

Contrary to the SD method, the minimal gradient (MG) method \cite{dai2003altermin} computes its stepsize by minimizing the gradient norm,
\begin{equation}\label{mg}
  \alpha_k^{MG}=\arg\min_{\alpha}~\|g(x_k-\alpha g_k)\|=\frac{g_{k} \tr Ag_{k}}{g_{k} \tr A^2g_{k}}.
  %=\frac{g_{k} \tr Ag_{k}}{g_{k} \tr A^2g_{k}}.
\end{equation}
It is widely accepted that the MG method can also perform poorly and  has similar asymptotic behavior as the SD method, i.e., it will asymptotically zigzag in a two-dimensional subspace.
In \cite{zhou2006gradient}, the authors provide some interesting analyses on $\alpha_k^{MG}$ for minimizing two-dimensional quadratics. However, rigorous asymptotic convergence results of the MG method
for minimizing general quadratic function are very limit in literature.

In order to avoid the zigzagging pattern, it is useful to determine the stepsize without using the exact stepsize because it would yield a  gradient perpendicular to the current one.
Barzilai and Borwein \cite{barzilai1988two} proposed the following two novel stepsizes:
\begin{equation}\label{sBB}
  \alpha_k^{BB1}=\frac{s_{k-1} \tr s_{k-1}}{s_{k-1} \tr y_{k-1}}~~\textrm{and}~~
  \alpha_k^{BB2}=\frac{s_{k-1} \tr y_{k-1}}{y_{k-1} \tr y_{k-1}},
\end{equation}
%\begin{equation*}
%  \alpha_k^{BB2}=\frac{s_{k-1} \tr y_{k-1}}{y_{k-1} \tr y_{k-1}},
%\end{equation*}
where $s_{k-1}=x_k-x_{k-1}$ and $y_{k-1}=g_k-g_{k-1}$. The BB method \eqref{sBB} performs quite well in practice, though it generates a nonmonotone sequence of objective values. Due to its simplicity and efficiency, the BB method has been widely studied \cite{dai2005asymptotic,dhl2018,dai2002r,fletcher2005barzilai,raydan1993barzilai} and extended to general problems and various applications, see \cite{birgin2000nonmonotone,huang2016smoothing,huang2015quadratic,jiang2013feasible,liu2011coordinated,raydan1997barzilai}.
%unconstrained problems \cite{raydan1997barzilai}, to constrained optimization problems \cite{birgin2000nonmonotone,huang2016smoothing}, and to various applications \cite{huang2015quadratic,jiang2013feasible,liu2011coordinated}.
Another line of research to break the zigzagging pattern and accelerate the convergence is occasionally applying short stepsizes that approximate $1/\lambda_n$
to eliminate the corresponding component of the gradient. One seminal work is due to Yuan \cite{yuan2006new,yuan2008step}, who derived the following stepsize:
\begin{equation}\label{syv}
  \alpha_k^{Y}=\frac{2}{\frac{1}{\alpha_{k-1}^{SD}}+\frac{1}{\alpha_{k}^{SD}}+
  \sqrt{\left(\frac{1}{\alpha_{k-1}^{SD}}-\frac{1}{\alpha_{k}^{SD}}\right)^2+
  \frac{4\|g_k\|^2}{(\alpha_{k-1}^{SD}\|g_{k-1}\|)^2}}}.
\end{equation}
Dai and Yuan \cite{dai2005analysis} further suggested a new gradient method with
\begin{equation}\label{sdy}
  \alpha_k^{DY}=\left\{
                  \begin{array}{ll}
                    \alpha_k^{SD}, & \hbox{if mod($k$,4)$<2$;} \\
                    \alpha_k^{Y}, & \hbox{otherwise.}
                  \end{array}
                \right.
\end{equation}
The DY method \eqref{sdy} is a monotone method and appears very competitive with the nonmonotone BB method. Recently, by employing the results in \cite{akaike1959successive,nocedal2002behavior}, De Asmundis et al. \cite{de2014efficient} show that the stepsize $\alpha_k^{Y}$ converges to $1/\lambda_n$ if the SD method is applied to problem \eqref{qudpro}. This spectral property is the key to break the zigzagging pattern.

In \cite{dai2006new2}, Dai and Yang developed the asymptotic optimal gradient (AOPT) method whose stepsize is given by
\begin{equation}\label{saopt}
  \alpha_k^{AOPT}=\frac{\|g_k\|}{\|Ag_k\|}.
\end{equation}
Unlike the DY method, the AOPT method only has one stepsize. In addition, they show that $\alpha_k^{AOPT}$ asymptotically converges to $\frac{2}{\lambda_1+\lambda_n}$, which is in some sense an optimal stepsize since it minimizes $\|I-\alpha A\|$ over $\alpha$ \cite{dai2006new2,elman1994inexact}. However, the AOPT method also asymptotically alternates between two directions. To accelerate the AOPT method, Huang et al. \cite{huang2019gradient} derived a new stepsize that converges to $1/\lambda_n$ during the AOPT iterates and further suggested a gradient method to exploit spectral properties of the stepsizes. For the latest developments of exploiting spectral properties to accelerate gradient methods, see \cite{de2014efficient,de2013spectral,di2018steplength,gonzaga2016steepest,huang2019gradient}.

In this paper, we present the analysis on the asymptotic behaviors of gradient methods and the techniques for breaking the zigzagging pattern. For a uniform analysis, we consider the following stepsize
\begin{equation}\label{glstep1}
  \alpha_k=\frac{g_k \tr \Psi(A) g_k}{g_k \tr \Psi(A)A g_k},
\end{equation}
where $\Psi$ is a real analytic function on $[\lambda_1,\lambda_n]$ and can be expressed by Laurent series
\begin{equation*}
  \Psi(z)=\sum_{k=-\infty}^\infty c_kz^k,~~c_k\in\mathbb{R},
\end{equation*}
such that $0<\sum_{k=-\infty}^\infty c_kz^k<+\infty$ for all $z\in[\lambda_1,\lambda_n]$. Apparently, $\alpha_k$ is a family of stepsizes that would give a family of gradient methods.
When $\Psi(A)=A^u$ for some nonnegative integer $u$, we get the following stepsize
\begin{equation}\label{glstep2}
  \alpha_k=\frac{g_k \tr A^{u} g_k}{g_k \tr A^{u+1} g_k}.
\end{equation}
The $\alpha_k^{SD}$ and $\alpha_k^{MG}$ simply correspond to the cases $u=0$ and $u=1$, respectively.

We will present theoretical analysis on the asymptotic convergence on the family of gradient methods whose stepsize can be written in the form \eqref{glstep1},
which  provides justifications for the zigzag behaviors of all these gradient methods including the SD and MG methods.
In particular, we show that each method in the family \eqref{glstep1} will asymptotically alternate between two directions associated with the two eigenvectors corresponding to $\lambda_1$ and $\lambda_n$. Moreover, we analyze the asymptotic behaviors of the objective value, gradient norm, and stepsize.
It is shown that, when $\Psi(A)\neq I$, the two sequences
$\Big\{\frac{\Delta_{2k+1}}{\Delta_{2k}} \Big\}$ and $\Big\{ \frac{\Delta_{2k+2}}{\Delta_{2k+1}} \Big\}$
may converge at different speeds, while the odd and even subsequences
$\Big\{\frac{\Delta_{2k+3}}{\Delta_{2k+1}}\Big\}$ and $\Big\{ \frac{\Delta_{2k+2}}{\Delta_{2k}}\Big\}$
converge at the same rate, where $\Delta_k = f(x_k) - f^* $.
%In addition, $\Big\{\frac{\|g_{2k+2}\|^2}{\|g_{2k}\|^2}\Big\}$ and $\Big\{\frac{\|g_{2k+3}\|^2}{\|g_{2k+1}\|^2}\Big\}$ also converge at the same rate.
Similar property is also possessed by the gradient norm sequence. In addition, we show each method in \eqref{glstep1} has the same worst asymptotic rate.

In order to accelerate the gradient methods \eqref{glstep1}, we investigate techniques for breaking the zigzagging pattern. We derive a new stepsize $\tilde{\alpha}_k$ based on finite termination for
minimizing two-dimensional strictly convex quadratic function. For the $n$-dimensional case,
we prove that $\tilde{\alpha}_k$ converges to $1/\lambda_n$ when gradient methods \eqref{glstep1} are applied to problem \eqref{qudpro}. Furthermore, based on this spectral property, we propose a periodic
gradient method, which, in a periodic mode, alternately uses the BB stepsize, stepsize \eqref{glstep1}
 and our new stepsize $\tilde{\alpha}_k$. Numerical comparisons of the proposed method with the BB \cite{barzilai1988two}, DY \cite{dai2005analysis}, ABBmin2 \cite{frassoldati2008new}, and SDC \cite{de2014efficient} methods show that the new gradient method is very efficient.
Our theoretical results also significantly improve and generalize those in \cite{akaike1959successive,nocedal2002behavior}, where only the SD method (i.e., $\Psi(A)=I$) is considered.
 We point out that \cite{pronzato2006asymptotic} does not analyze the asymptotic behaviors of the objective value, gradient norm, and stepsize, though \eqref{glstep1} is similar to the $P$-gradient methods in \cite{pronzato2006asymptotic}. Moreover, we  develop techniques for accelerating these zigzag methods with simpler analysis. Notice that $\alpha_k^{AOPT}$ can not be written in the form \eqref{glstep1}. Thus, our results are not applicable to the AOPT method. On the other hand, the analysis of the AOPT method presented in
\cite{dai2006new2}  can not be applied directly to the family of methods \eqref{glstep1}.

%\cite{akaike1959successive,forsythe1968asymptotic,gonzaga2016steepest,nocedal2002behavior} for more detail.

The paper is organized as follows. In Section \ref{secasmg}, we analyze the asymptotic behaviors of the family
of gradient methods \eqref{glstep1}.
In Section \ref{secnewstep}, we accelerate the gradient methods \eqref{glstep1} by developing techniques to
break its zigzagging pattern and propose a new periodic gradient method.
Numerical experiments are presented in Section \ref{secnum}.
 Finally, some conclusions and discussions are made in Section \ref{seccls}.

\section{Asymptotic behavior of the family \eqref{glstep1}}\label{secasmg}
%In this section, inspired \cite{akaike1959successive,dai2006new2},
In this section, we present a uniform analysis on the asymptotic behavior of the family of gradient methods \eqref{glstep1} for general $n$-dimensional strictly convex quadratics.

%We mention that the family \eqref{glstep1} is monotone when the minimum order of the series of $\Psi(A)$ is positive, and our analysis applies to all monotone gradient methods whose stepsize can be written in the form \eqref{glstep1}, such as the SD and MG methods.

Let $\{\lambda_1,\lambda_2,\cdots,\lambda_n\}$ be the eigenvalues of $A$, and $\{\xi_1,\xi_2,\ldots,\xi_n\}$ be the associated orthonormal eigenvectors. Noting that the gradient method is invariant under translations and rotations when applying to a quadratic function. For theoretical analysis, we can assume without loss of generality that
\begin{equation}\label{formA}
  A=\textrm{diag}\{\lambda_1,\lambda_2,\cdots,\lambda_n\},~~0<\lambda_1<\lambda_2<\cdots<\lambda_n.
\end{equation}

Denoting the components of $g_k$ along the eigenvectors $\xi_i$ by $\mu_k^{(i)}$, $i=1,\ldots,n$, i.e.,
\begin{equation}\label{gmu}
  g_k=\sum_{i=1}^n\mu_k^{(i)}\xi_i.
\end{equation}
The above decomposition of gradient $g_k$ together with the update rule \eqref{eqitr} gives that
\begin{equation}\label{gupeq}
  g_{k+1}=g_k-\alpha_kAg_k=\prod_{j=1}^k(I-\alpha_jA)g_0=\sum_{i=1}^n\mu_{k+1}^{(i)}\xi_i,
\end{equation}
where
\begin{equation}\label{upmuk}
  \mu_{k+1}^{(i)}=(1-\alpha_k\lambda_i)\mu_k^{(i)}=\mu_0^{(i)}\prod_{j=1}^k(1-\alpha_j\lambda_i).
\end{equation}
%This relation implies that the closer $\alpha_k$ to $\frac{1}{\lambda_i}$, the smaller $|\mu_i^{k+1}|$ would be. In addition, if $\mu_i^{k}=0$, the corresponding component will vanish at all subsequent iterations.{upmuk}
Defining the vector $q_k=\left(q_k^{(i)}\right)$ with
\begin{equation}\label{defpk}
  q_k^{(i)}=\frac{(\mu_k^{(i)})^2}{\|\mu_k\|^2}
\end{equation}
and
\begin{equation}\label{defgama}
  \gamma_k=\frac{1}{\alpha_k}=\frac{g_k \tr \Psi(A)A g_k}{g_k \tr \Psi(A) g_k}
  =\frac{\sum_{i=1}^n\Psi(\lambda_i)\lambda_iq_k^{(i)}}{\sum_{i=1}^n\Psi(\lambda_i)q_k^{(i)}},
\end{equation}
we can have from \eqref{upmuk}, \eqref{defpk} and \eqref{defgama} that
\begin{equation}\label{pk1}
  q_{k+1}^{(i)}=\frac{(\lambda_i-\gamma_k)^2q_k^{(i)}}
  {\sum_{i=1}^n(\lambda_i-\gamma_k)^2q_k^{(i)}}.
\end{equation}
In addition, by the definition of $q_k$, we know that $q_k^{(i)}\geq0$ for all $i$ and
\[
  \sum_{i=1}^nq_k^{(i)}=1,~~\forall~~k\geq1.
\]

%We further assume that the starting point $x_1$ is such that
%\begin{equation*}
%  g_0 \tr \xi_1\neq0,~~g_0 \tr \xi_n\neq0.
%\end{equation*}

Before establishing the asymptotic convergence of the family of gradient methods \eqref{glstep1},
we first give some lemmas on the properties of the sequence $\{q_k\}$.

\begin{lemma}\label{lm1}
  Suppose $p\in\mathbb{R}^n$ satisfies (i) $p^{(i)}\geq0$ for all $i=1,2,\ldots,n$; (ii) there exist at least two $i's$ with $p^{(i)}>0$; and (iii) $\sum_{i=1}^np^{(i)}=1$.
  Define $T: \mathbb{R}^n\rightarrow\mathbb{R}$ be the following transformation:
  \begin{equation}\label{tranT}
    (Tp)^{(i)}=\frac{(\lambda_i-\gamma(p))^2p^{(i)}}
    {\sum_{i=1}^n(\lambda_i-\gamma(p))^2p^{(i)}},
\end{equation}
  where
  \begin{equation}\label{defgama2}
    \gamma(p)=\frac{\sum_{i=1}^n\Psi(\lambda_i)\lambda_ip^{(i)}}{\sum_{i=1}^n\Psi(\lambda_i)p^{(i)}}.
  \end{equation}
  Then we have
  \begin{equation}\label{thtinq}
    \Theta(Tp)\geq\Theta(p),
  \end{equation}
  where
    \begin{equation}\label{deftht}
    \Theta(p)=\frac{\sum_{i=1}^n\Psi(\lambda_i)(\lambda_i-\gamma(p))^2p^{(i)}}
    {\sum_{i=1}^n\Psi(\lambda_i)p^{(i)}}.
  \end{equation}
  In addition, \eqref{thtinq} holds with equality if and only if there are two indices, say $i_1$ and $i_2$, such that $p^{(i)}=0$ for all $i\notin\{i_1,i_2\}$ and
  \begin{equation}\label{sum2eig}
\gamma(Tp)+\gamma(p)=\lambda_{i_1}+\lambda_{i_2}.
\end{equation}
\end{lemma}
\begin{proof}
It follows from the definition of $Tp$ that
    \begin{align}\label{thtp}
    \Theta(Tp)&=\frac{\sum_{i=1}^n\Psi(\lambda_i)(\lambda_i-\gamma(Tp))^2(Tp)^{(i)}}
    {\sum_{i=1}^n\Psi(\lambda_i)(Tp)^{(i)}}
    \nonumber\\
    &=\frac{\sum_{i=1}^n\Psi(\lambda_i)(\lambda_i-\gamma(Tp))^2(\lambda_i-\gamma(p))^2p^{(i)}}
    {\sum_{i=1}^n\Psi(\lambda_i)(\lambda_i-\gamma(p))^2p^{(i)}}.
\end{align}
Let us define two vectors $w = (w_i) \in \mathbb{R}^n$ and $z = (z_i) \in \mathbb{R}^n$ by
\[
  w_i=\sqrt{\Psi(\lambda_i)}(\lambda_i-\gamma(Tp))(\lambda_i-\gamma(p))\sqrt{p^{(i)}}
\]
and
\[
  z_i=\sqrt{\Psi(\lambda_i)}\sqrt{p^{(i)}}.
\]
Then, we have from the Cauchy-Schwarz inequality that
\begin{align}\label{thtpinq}
  \|w\|^2\|z\|^2&=\left(\sum_{i=1}^n\Psi(\lambda_i)(\lambda_i-\gamma(Tp))^2(\lambda_i-\gamma(p))^2p^{(i)}\right)
  \left(\sum_{i=1}^n\Psi(\lambda_i)p^{(i)}\right)
  \nonumber\\
  &\geq (w \tr z)^2=\left(\sum_{i=1}^n\Psi(\lambda_i)(\lambda_i-\gamma(Tp))(\lambda_i-\gamma(p))p^{(i)}\right)^2.
\end{align}
Using the definition of $\gamma(p)$, we can obtain that
\begin{eqnarray}\label{keyeq}
 && \sum_{i=1}^n\Psi(\lambda_i)(\lambda_i-\gamma(Tp))(\lambda_i-\gamma(p))p^{(i)}
  -\sum_{i=1}^n\Psi(\lambda_i)(\lambda_i-\gamma(p))^2p^{(i)} \nonumber \\
 & = &(\gamma(p)-\gamma(Tp))\sum_{i=1}^n\Psi(\lambda_i)(\lambda_i-\gamma(p))p^{(i)}
  =0,
\end{eqnarray}
which together with \eqref{thtpinq} gives
\begin{eqnarray}\label{thtpinq2}
 && \left(\sum_{i=1}^n\Psi(\lambda_i)(\lambda_i-\gamma(Tp))^2(\lambda_i-\gamma(p))^2p^{(i)}\right)
  \left(\sum_{i=1}^n\Psi(\lambda_i)p^{(i)}\right)
  \nonumber\\
 & \geq & \left(\sum_{i=1}^n\Psi(\lambda_i)(\lambda_i-\gamma(p))^2p^{(i)}\right)^2.
\end{eqnarray}
Then, the inequality \eqref{thtinq} follows immediately.

The equality in \eqref{thtpinq} holds if and only if
\begin{equation}\label{equv2}
\sqrt{\Psi(\lambda_i)}(\lambda_i-\gamma(Tp))(\lambda_i-\gamma(p))\sqrt{p^{(i)}}
  =C\sqrt{\Psi(\lambda_i)}\sqrt{p^{(i)}},~~i=1,\ldots,n
\end{equation}
for some nonzero scalar $C$. Clearly, \eqref{equv2} holds when $p^{(i)}=0$.
Suppose that there exist two indices $i_1$ and $i_2$ such that $p^{(i_1)},p^{(i_2)}>0$.
It follows from \eqref{equv2} that
\[
(\lambda_{i_1}-\gamma(Tp))(\lambda_{i_1}-\gamma(p))
  =(\lambda_{i_2}-\gamma(Tp))(\lambda_{i_2}-\gamma(p)).
\]
%i.e.,
%\begin{equation*}
%\frac{\lambda_{i_1}-\gamma(Tp)}
%{\lambda_{i_2}-\gamma(Tp)}
%=\frac{\lambda_{i_2}-\gamma(p)}{\lambda_{i_1}-\gamma(p)}.
%\end{equation*}
So, by the assumption \eqref{formA}, we have
\[
\lambda_{i_1}+\lambda_{i_2}=\gamma(Tp)+\gamma(p),
\]
which again with assumption \eqref{formA} imply that
 \eqref{equv2} holds  if and only if $p$ has only two nonzero components
 and \eqref{sum2eig} holds.
\end{proof}

\begin{lemma}\label{tpeqp}
Let $p_*\in\mathbb{R}^n$ satisfy the conditions of Lemma \ref{lm1} and $T$ be the transformation \eqref{tranT}.
If $p_*$ has only two nonzero components $p_*^{(i_1)}$ and $p_*^{(i_2)}$, we have
    \begin{equation}\label{tps1}
      (Tp_*)^{(i_1)}=\frac{\Psi^2(\lambda_{i_2})p_*^{(i_2)}}
      {\Psi^2(\lambda_{i_1})p_*^{(i_1)}+\Psi^2(\lambda_{i_2})p_*^{(i_2)}},
  \end{equation}
      \begin{equation}\label{tps2}
   (Tp_*)^{(i_2)}=\frac{\Psi^2(\lambda_{i_1})p_*^{(i_1)}}
      {\Psi^2(\lambda_{i_1})p_*^{(i_1)}+\Psi^2(\lambda_{i_2})p_*^{(i_2)}},
  \end{equation}
  \begin{equation}\label{tsqstar}
    (T^2p_*)^{(i_1)}=p_*^{(i_1)},~~(T^2p_*)^{(i_2)}=p_*^{(i_2)},
  \end{equation}
    and
  \begin{equation}\label{gmaptp}
    \gamma(p_*)+\gamma(Tp_*)=\lambda_{i_1}+\lambda_{i_2},
  \end{equation}
where the function $\gamma$ is defined in \eqref{defgama2}.
Moreover, $p_*=Tp_*$ if and only if
  \begin{equation}\label{ps1}
    p_*^{(i_1)}=\frac{\Psi(\lambda_{i_2})}{\Psi(\lambda_{i_1})+\Psi(\lambda_{i_2})} \quad \mbox{and} \quad
  p_*^{(i_2)}=\frac{\Psi(\lambda_{i_1})}{\Psi(\lambda_{i_1})+\Psi(\lambda_{i_2})}.
  \end{equation}
\end{lemma}
\begin{proof}
  %The conclusions follow by direct computation.
  By the definition of $\gamma(p)$, we have
  \begin{equation}\label{gmsr1}
    \gamma(p_*)=\frac{\Psi(\lambda_{i_1})\lambda_{i_1}p_*^{(i_1)}+\Psi(\lambda_{i_2})\lambda_{i_2}p_*^{(i_2)}}
    {\Psi(\lambda_{i_1})p_*^{(i_1)}+\Psi(\lambda_{i_2})p_*^{(i_2)}},
  \end{equation}
  which indicates that
  \[
    \lambda_{i_1}-\gamma(p_*)=\frac{\Psi(\lambda_{i_2})p_*^{(i_2)}(\lambda_{i_1}-\lambda_{i_2})}
    {\Psi(\lambda_{i_1})p_*^{(i_1)}+\Psi(\lambda_{i_2})p_*^{(i_2)}}, \quad
\lambda_{i_2}-\gamma(p_*)=\frac{\Psi(\lambda_{i_1})p_*^{(i_1)}(\lambda_{i_2}-\lambda_{i_1})}
    {\Psi(\lambda_{i_1})p_*^{(i_1)}+\Psi(\lambda_{i_2})p_*^{(i_2)}}.
  \]
 % and
 %   \begin{equation*}
 %   \lambda_{i_2}-\gamma(p_*)=\frac{\Psi(\lambda_{i_1})p_*^{(i_1)}(\lambda_{i_2}-\lambda_{i_1})}
 %   {\Psi(\lambda_{i_1})p_*^{(i_1)}+\Psi(\lambda_{i_2})p_*^{(i_2)}}.
 % \end{equation*}
Then, it follows from the definition of transformation $T$ that
  \begin{align*}
    (Tp_*)^{(i_1)}&=\frac{(\Psi(\lambda_{i_2})p_*^{(i_2)})^2p_*^{(i_1)}}
    {(\Psi(\lambda_{i_2})p_*^{(i_2)})^2p_*^{(i_1)}+(\Psi(\lambda_{i_1})p_*^{(i_1)})^2p_*^{(i_2)}}
    \\
    &=\frac{\Psi^2(\lambda_{i_2})p_*^{(i_2)}}
      {\Psi^2(\lambda_{i_1})p_*^{(i_1)}+\Psi^2(\lambda_{i_2})p_*^{(i_2)}}.
  \end{align*}
  This gives \eqref{tps1}. \eqref{tps2} can be proved similarly.
By \eqref{tps1} and \eqref{tps2}, we have
    \begin{align*}
    (T^2p_*)^{(i_1)}
    &=\frac{\Psi^2(\lambda_{i_2})(Tp_*)^{(i_2)}}
      {\Psi^2(\lambda_{i_1})(Tp_*)^{(i_1)}+\Psi^2(\lambda_{i_2})(Tp_*)^{(i_2)}}\\
    &=\frac{\Psi^2(\lambda_{i_1})\Psi^2(\lambda_{i_2})p_*^{(i_1)}}
    {\Psi^2(\lambda_{i_1})\Psi^2(\lambda_{i_2})p_*^{(i_2)}+
    \Psi^2(\lambda_{i_1})\Psi^2(\lambda_{i_2})p_*^{(i_1)}}\\
    & =\frac{p_*^{(i_1)}}
      {p_*^{(i_1)}+p_*^{(i_2)}}=p_*^{(i_1)}.
  \end{align*}
  $(T^2p_*)^{(i_2)}$ follows similarly. This proves \eqref{tsqstar}.

Again by \eqref{tps1}, \eqref{tps2} and  the definition of function $\gamma$ in \eqref{defgama2}, we have
    \begin{equation}\label{gmsrtp}
    \gamma(Tp_*)=\frac{\lambda_{i_1}\Psi(\lambda_{i_2})p_*^{(i_2)}+
    \lambda_{i_2}\Psi(\lambda_{i_1})p_*^{(i_1)}}
    {\Psi(\lambda_{i_1})p_*^{(i_1)}+\Psi(\lambda_{i_2})p_*^{(i_2)}}.
  \end{equation}
Then, the equality \eqref{gmaptp} follows from \eqref{gmsr1} and \eqref{gmsrtp}.
 For \eqref{ps1}, let
  \[
    p_*^{(i_1)}=\frac{\Psi^2(\lambda_{i_2})p_*^{(i_2)}}
      {\Psi^2(\lambda_{i_1})p_*^{(i_1)}+\Psi^2(\lambda_{i_2})p_*^{(i_2)}}.
  \]
  Rearranging terms and using $p_*^{(i_1)}+p_*^{(i_2)}=1$, we have
  \[
   \Psi^2(\lambda_{i_1})(p_*^{(i_1)})^2=\Psi^2(\lambda_{i_2})(p_*^{(i_2)})^2,
  \]
  which implies that
  \[
    \Psi(\lambda_{i_1})p_*^{(i_1)}=\Psi(\lambda_{i_2})p_*^{(i_2)}.
  \]
  This together with the fact $p_*^{(i_1)}+p_*^{(i_2)}=1$ yields \eqref{ps1}.
\end{proof}

%For any $p\in\mathbb{R}^n$, denote $T^0_p=p$. We have the following lemma.
\begin{lemma}\label{lm3}
Let $p\in\mathbb{R}^n$ satisfy the conditions of Lemma \ref{lm1} and $T$ be the transformation  \eqref{tranT}. Then, there exists a $p_*$ satisfying
\begin{equation}\label{t2k}
  \lim_{k\rightarrow\infty}T^{2k}p=p_*~~and~~
  \lim_{k\rightarrow\infty}T^{2k+1}p=Tp_*,
  \end{equation}
where $p_*$ and $Tp_*$ have only two nonzero components satisfying
\begin{equation}\label{limp3}
  p_*^{(i_1)}+p_*^{(i_2)}=1, ~~p_*^{(i)}=0,~~i\neq i_1,i_2,
\end{equation}
\begin{equation}\label{limp4}
  (Tp_*)^{(i_1)}+(Tp_*)^{(i_2)}=1,~~(Tp_*)^{(i)}=0,~~i\neq i_1,i_2,
\end{equation}
for some $i_1,i_2\in\{1,\ldots,n\}$. Hence, \eqref{tps1}, \eqref{tps2}, \eqref{tsqstar} and \eqref{gmaptp} hold.
\end{lemma}
\begin{proof}
Let $p_0=T^0p=p$ and $p_{k}=T{p_{k-1}}=T^{k}{p_0}$. Obviously, for all $k\geq0$, $p_{k}$ satisfies (i) and (iii) of Lemma \ref{lm1}. Let $i_{\min} = \min\{i \in \mathcal{N}: p_0^{(i)}>0\}$
and $i_{\max} = \max \{i \in \mathcal{N}: p_0^{(i)}>0\}$, where $\mathcal{N} = \{1, \ldots, n\}$.
%Suppose that $i_{\min}$ and $i_{\max}$ are the minimal and maximal superscripts such that $p_0^{(i)}>0$.
From the definition of $\gamma$, we know $\lambda_{i_{\min}}< \gamma(p)<\lambda_{i_{\max}}$.
Thus, by the definition of $T$, we have $p_1^{({i_{\min}})}>0$ and $p_1^{({i_{\max}})}>0$. Then, by induction, for all $k\geq0$, $p_{k}$ satisfies (ii) of Lemma \ref{lm1}.
So, by Lemma \ref{lm1},  $\{\Theta(p_k)\}$ is a monotonically increasing sequence.  Since  $\lambda_1\leq \gamma(p)\leq\lambda_n$,
we have $(\lambda_i-\gamma(p))^2\leq(\lambda_n-\lambda_1)^2$.
 Hence, we have from the definition of $\Theta$ that
%  \Theta(p_k)\leq\frac{(\lambda_n-\lambda_1)^2\max_i\{\Psi(\lambda_i)\}}{\min_i\{\Psi(\lambda_i)\}}.
$\Theta(p_k)\leq (\lambda_n-\lambda_1)^2 $.
%$\Theta(p_k)\leq\frac{\lambda_n^{u}(\lambda_n-\lambda_1)^2}{\lambda_1^{u}}$.
Thus, $\{\Theta(p_k)\}$ is convergent. Let $\Theta_*=\lim_{k\rightarrow\infty}\Theta(p_k) >0$.

Denote the set of all limit points of $\{p_k\}$ by $P_*$ with cardinality $|P_*|$. Since $\{p_k\}$ is bounded, $|P_*|\geq1$.
For any subsequence $\{p_{k_j}\}$ converging to some $p_*\in P_*$, we have
\[
    \lim_{j\rightarrow\infty}\Theta(p_{k_j})=\Theta(p_*) \quad \mbox{and} \quad   \lim_{j\rightarrow\infty}\Theta(Tp_{k_j})=\Theta(Tp_*),
 \]
 by the continuity of $\Theta$ and $T$.
Notice $p_{k_j+1}=Tp_{k_j}$, we have $ \Theta_*=\Theta(p_*)=\Theta(Tp_*)$.

Since $p_{k}$ satisfies (i)-(iii) of Lemma \ref{lm1} for all $k\geq0$, $p_*$ must satisfy (i) and (iii). If $p_*$ has only one positive component, we have $\Theta(p_*)=0$
which contradicts $\Theta(p_*) = \Theta_* > 0 $.
Hence, by Lemma \ref{lm1}, Lemma \ref{tpeqp} and $\Theta(p_*)=\Theta(Tp_*)$, $p_*$ has only two nonzero components,
say $p_*^{(i_1)}$ and $p_*^{(i_2)}$, and their values are uniquely determined by the indices $i_1$, $i_2$ and the eigenvalues $\lambda_{i_1}$ and $\lambda_{i_2}$.
This implies $|P_*| < \infty$. Furthermore, by Lemma \ref{tpeqp}, for any $p_*\in P_*$, $Tp_*$ is given by \eqref{tps1} and \eqref{tps2}, and $Tp_*\in P_*$.

We now show that $|P_*|\leq2$ by way of contradiction.  Suppose $|P_*| \ge 3 $. For any $p_*\in P_*$ and $Tp_*\in P_*$, denote $\delta_1$ and $\delta_2$ to be the distance from $p_*$ to $P_*\setminus \{p_*\}$ and from $Tp_*$ to $P_*\setminus \{Tp_*\}$, respectively. Since $ 3 \le |P_*| < \infty$, we have $\delta_1>0$, $\delta_2>0$ and there exists an infinite subsequence $\{p_{k_j}\}$ such that
  \[
    p_{k_j}\rightarrow p_*, \quad \mbox{and} \quad
    p_{k_j+1}=Tp_{k_j}\rightarrow Tp_*,
  \]
  but $p_{k_j+2}\notin\mathcal{B}\left(p_*,\frac{1}{2}\delta\right)\cup
  \mathcal{B}\left(Tp_*,\frac{1}{2}\delta\right)$, where $\delta=\min\{\delta_1,\delta_2\}$ and $\mathcal{B}(p_*,r)=\{p: \|p-p_*\|\leq r\}$.
However, by \eqref{tsqstar} we have $T^2p_*=p_*$. Hence, by continuity of $T$,
  \[
    \lim_{j\rightarrow\infty}p_{k_j+2}=\lim_{j\rightarrow\infty}Tp_{k_j+1}=\lim_{j\rightarrow\infty}T^2p_{k_j}= p_*,
\]
which contradicts the choice of $p_{k_j+2}\notin\mathcal{B}\left(p_*,\frac{1}{2}\delta\right)$.
Thus, $\{p_k\}$ has at most two limit points $p_*$ and $Tp_*$, and both have only two nonzero components.

Now, we assume that $p_*$ is a limit point of $\{p_{2k}\}$. Since $T^2p_*=p_*$, all subsequences of $\{p_{2k}\}$ have the same limit point, i.e.,
$p_{2k}=T^{2k}p\rightarrow p_*$. Similarly, we have $T^{2k+1}p\rightarrow Tp_*$. Then, \eqref{limp3} and \eqref{limp4} follow directly from the analysis.
\end{proof}

%Notice that,, we have $q_{k+1}=Tp$ when $p=q_k$.
%Clearly, if $p=q_k$, we have by the definitions of $T$, $\gamma$, and $q_k$ that $q_{k+1}=Tp$. Hence, if $\{q_k\}$ satisfies the conditions in Lemma \ref{lm1}, the conclusions in Lemma \ref{lm3} hold.

Based on the above analysis, we can show that each gradient method in \eqref{glstep1} will asymptotically reduces its search
in a two-dimensional subspace spanned by the two eigenvectors $\xi_1$ and $\xi_n$.
\begin{theorem}\label{th1}
 Assume that the starting point $x_0$ has the property that
\begin{equation}\label{assp1}
  g_0 \tr \xi_1\neq0~~and~~g_0 \tr \xi_n\neq0.
\end{equation}
Let $\{x_k\}$ be the iterations generated by applying a method in \eqref{glstep1} to solve problem \eqref{qudpro}. Then
 \begin{equation}\label{mu2k}
  \lim_{k\rightarrow\infty}\frac{(\mu_{2k}^{(i)})^2}{\sum_{j=1}^n(\mu_{2k}^{(j)})^2}=
  \left\{
    \begin{array}{ll}
     \displaystyle\frac{1}{1+c^2}, & \hbox{if $i=1$,} \\
      0, & \hbox{if $i=2,\ldots,n-1$,} \\
      \displaystyle\frac{c^2}{1+c^2}, & \hbox{if $i=n$,}
    \end{array}
  \right.
  \end{equation}
and
\begin{equation}\label{mu2k1}
  \lim_{k\rightarrow\infty}\frac{(\mu_{2k+1}^{(i)})^2}{\sum_{j=1}^n(\mu_{2k+1}^{(j)})^2}=
  \left\{
       \begin{array}{ll}
        \displaystyle\frac{c^2\Psi^2(\lambda_{n})}{\Psi^2(\lambda_{1})+c^2\Psi^2(\lambda_{n})}, & \hbox{if $i=1$,} \\
         0, & \hbox{if $i=2,\ldots,n-1$,} \\
         \displaystyle\frac{\Psi^2(\lambda_{1})}{\Psi^2(\lambda_{1})+c^2\Psi^2(\lambda_{n})}, & \hbox{if $i=n$,}
       \end{array}
     \right.
   \end{equation}
where $c$ is a nonzero constant.
\end{theorem}
\begin{proof}
By the assumption \eqref{assp1}, we know that $q_0$ satisfies (i)-(iii) of Lemma \ref{lm1}.
Notice that $q_k = T^k q_0$. Then, by Lemma \ref{lm3}, there exists a $p_*$ such that the sequences $\{q_{2k}\}$ and $\{q_{2k+1}\}$ converge to $p_*$ and $Tp_*$,
respectively, which have only two nonzero components satisfying \eqref{limp3}, \eqref{limp4}  for some $i_1,i_2\in\{1,\ldots,n\}$, and \eqref{tsqstar} holds.
 Hence, if $1 \le i_1<i_2<n$, we have
\begin{equation}\label{p2klim}
  \lim_{k\rightarrow\infty}q_{2k}^{(n)}=0, \qquad  \lim_{k\rightarrow\infty}\frac{q_{2k}^{(i_2)}}{q_{2k+2}^{(i_2)}}=1,
\end{equation}
and
\[
\lim_{k\rightarrow\infty}(\gamma(q_{2k})+\gamma(q_{2k+1}))= \gamma(p_*)+\gamma(Tp_*) = \lambda_{i_1}+\lambda_{i_2}.
\]
In addition, since $q_0^{(1)} > 0 $ and $q_0^{(n)} > 0$ by \eqref{assp1},
 we can see from  the proof of Lemma \ref{lm3} that  $q_k^{(1)}>0$, $q_k^{(n)}>0$ for all $k\geq0$. Thus, we have
\begin{align}
  \lim_{k\rightarrow\infty}\frac{q_{2k+2}^{(n)}}{q_{2k}^{(n)}}
&  =\lim_{k\rightarrow\infty}\frac{q_{2k+2}^{(n)}}{q_{2k}^{(n)}}\frac{q_{2k}^{(i_2)}}{q_{2k+2}^{(i_2)}}
 =\lim_{k\rightarrow\infty}\frac{(\lambda_n-\gamma(q_{2k+1}))^2(\lambda_n-\gamma(q_{2k}))^2}
{(\lambda_{i_2}-\gamma(q_{2k+1}))^2(\lambda_{i_2}-\gamma(q_{2k}))^2} \nonumber \\
&=\lim_{k\rightarrow\infty} \left(\frac{\lambda_n^2-(\gamma(q_{2k})+\gamma(q_{2k+1}))\lambda_n+\gamma(q_{2k})\gamma(q_{2k+1})}
{\lambda_{i_2}^2-(\gamma(q_{2k})+\gamma(q_{2k+1}))\lambda_{i_2}+\gamma(q_{2k})\gamma(q_{2k+1})}\right)^2 \nonumber \\
&=\left(\frac{\lambda_n^2-(\lambda_{i_1}+\lambda_{i_2})\lambda_n+\tilde{\gamma}}
{\lambda_{i_2}^2-(\lambda_{i_1}+\lambda_{i_2})\lambda_{i_2}+\gamma(p_*)\gamma(Tp_*)}\right)^2 \nonumber\\
&=\left(1+\frac{(\lambda_n-\lambda_{i_1})(\lambda_n-\lambda_{i_2})}
{\lambda_{i_2}^2-(\lambda_{i_1}+\lambda_{i_2})\lambda_{i_2}+\gamma(p_*)\gamma(Tp_*)}\right)^2 = : \rho. \label{huang111}
\end{align}
Since $\lambda_{i_1} < \gamma(p_*) <  \lambda_{i_2}$ and $\lambda_{i_1} < \gamma(Tp_*) <  \lambda_{i_2}$, we have
\begin{eqnarray*}
  \lambda_{i_2}^2-(\lambda_{i_1}+\lambda_{i_2})\lambda_{i_2}+\gamma(p_*)\gamma(Tp_*)
&= &\lambda_{i_2}^2-( \gamma(p_*)+\gamma(Tp_*) )\lambda_{i_2}+\gamma(p_*)\gamma(Tp_*)\\
&= & (\lambda_{i_2} - \gamma(p_*)) (\lambda_{i_2} - \gamma(Tp_*)) > 0.
\end{eqnarray*}
Hence, it follows from \eqref{huang111} that $\rho > 1$. So, $q_{2k}^{(n)}\rightarrow+\infty$, which contradicts \eqref{p2klim}.
Then, we must have $i_2=n$. In a similar way, we can show that $i_1=1$.
Finally, the equalities in \eqref{mu2k} and \eqref{mu2k1} follow directly from Lemma \ref{tpeqp}.
\end{proof}

In the following, we refer $c$ as the same constant in Theorem \ref{th1}.
By Theorem \ref{th1}  we can directly obtain the asymptotic behavior of the stepsize.
\begin{corollary}\label{cor1}
Under the conditions of Theorem \ref{th1}, we have
\begin{equation}\label{salp1}
  \lim_{k\rightarrow\infty}\alpha_{2k}=
  \frac{\Psi(\lambda_{1})+c^2\Psi(\lambda_{n})}{\lambda_{1}(\Psi(\lambda_{1})+c^2\kappa\Psi(\lambda_{n}))}
\end{equation}
and
\begin{equation}\label{salp2}
  \lim_{k\rightarrow\infty}\alpha_{2k+1}=
  \frac{\Psi(\lambda_{1})+c^2\Psi(\lambda_{n})}{\lambda_{1}(\kappa\Psi(\lambda_{1})+c^2\Psi(\lambda_{n}))},
\end{equation}
where $\alpha_k$ is defined in \eqref{glstep1} and $\kappa = \lambda_n/\lambda_1$ is the condition number of $A$. Moreover,
  \begin{equation}\label{sum2reps}
    \lim_{k\rightarrow\infty}\left(\frac{1}{\alpha_{2k}}+\frac{1}{\alpha_{2k+1}}\right)
    =\lambda_{1}+\lambda_{n}.
  \end{equation}
\end{corollary}

The next corollary interprets the constant $c$. A special result for the case $\Psi(A)=I$ (i.e., the SD method) can be found in Lemma 3.4 of \cite{nocedal2002behavior}.
\begin{corollary}\label{cvalue}
Under the conditions of Theorem \ref{th1}, we have
  \begin{equation}\label{eqcvalue}
  c=\lim_{k\rightarrow\infty}\frac{\mu_{2k}^{(n)}}{\mu_{2k}^{(1)}}
  =-\frac{\Psi(\lambda_{1})}{\Psi(\lambda_{n})}\lim_{k\rightarrow\infty}\frac{\mu_{2k+1}^{(1)}}{\mu_{2k+1}^{(n)}}.
\end{equation}
%if \eqref{mu2k1} holds, the constant $c$ satisfies
%  \begin{equation}\label{eqcvalue2}
%  c=\frac{1}{\kappa^u}\lim_{k\rightarrow\infty}\frac{\mu_{2k}^{(1)}}{\mu_{2k}^{(n)}}
%  =-\lim_{k\rightarrow\infty}\frac{\mu_{2k+1}^{(n)}}{\mu_{2k+1}^{(1)}}.
%\end{equation}
%In addition, $c$ is uniquely determined by the starting point $x_0$.
%and by the eigenvalues and the eigenvectors of $A$.
\end{corollary}
\begin{proof}
It follows from Theorem \ref{th1} that
  \begin{equation}\label{csqr}
    \lim_{k\rightarrow\infty}\frac{(\mu_{2k}^{(n)})^2}{(\mu_{2k}^{(1)})^2}
    =\frac{\Psi^2(\lambda_{1})}{\Psi^2(\lambda_{n})}
    \lim_{k\rightarrow\infty}\frac{(\mu_{2k+1}^{(1)})^2}{(\mu_{2k+1}^{(n)})^2}=c^2.
  \end{equation}
Note that $1/\lambda_n<\alpha_k<1/\lambda_1$ by the assumption \eqref{assp1}. And we have by \eqref{upmuk} that
  \begin{equation*}
    \mu_{2k+2}^{(1)}= \prod_{\ell=1}^2 (1-\alpha_{2k+\ell}\lambda_1)\mu_{2k}^{(1)} \quad \mbox{and} \quad
\mu_{2k+2}^{(n)}= \prod_{\ell=1}^2  (1-\alpha_{2k+\ell}\lambda_n)\mu_{2k}^{(n)}.
  \end{equation*}
 Thus, the sequence $\Big\{\frac{\mu_{2k}^{(n)}}{\mu_{2k}^{(1)}}\Big\}$, and similarly for $\Big\{\frac{\mu_{2k+1}^{(1)}}{\mu_{2k+1}^{(n)}}\Big\}$, do not change its sign.
Hence, without loss of generality, we can assume by \eqref{csqr} that
 \begin{equation}\label{eqcvalue2}
  %c=\lim_{k\rightarrow\infty}\frac{\mu_{2k}^{(n)}}{\mu_{2k}^{(1)}}.
  c=\lim_{k\rightarrow\infty} \mu_{2k}^{(n)} / \mu_{2k}^{(1)}.
\end{equation}
Then, by \eqref{upmuk}, \eqref{salp1} and \eqref{eqcvalue2},  we have
\[
    \lim_{k\rightarrow\infty}\frac{\mu_{2k+1}^{(1)}}{\mu_{2k+1}^{(n)}}
 =\lim_{k\rightarrow\infty}\frac{\mu_{2k}^{(1)}(1-\alpha_{2k}\lambda_1)}{\mu_{2k}^{(n)}(1-\alpha_{2k}\lambda_n)}
=-c\frac{\Psi(\lambda_{n})}{\Psi(\lambda_{1})},
\]
which gives \eqref{eqcvalue}.
\end{proof}

We have the following results on the asymptotic convergence of the function value.
\begin{theorem}\label{thconfk}
Under the conditions of Theorem \ref{th1}, we have
  \begin{equation}\label{f2k2f1k}
     \lim_{k\rightarrow\infty}\frac{f(x_{2k+1})-f^*}{f(x_{2k})-f^*}
     =R_f^1
\quad \mbox{and} \quad    \lim_{k\rightarrow\infty}\frac{f(x_{2k+2})-f^*}{f(x_{2k+1})-f^*}
     =R_f^2,
  \end{equation}
where
\begin{equation}\label{ratefk1}
  R_f^1=\frac{c^2(\kappa-1)^2(\Psi^2(\lambda_{1})+c^2\kappa\Psi^2(\lambda_{n}))}
  {(\Psi(\lambda_{1})+c^2\kappa\Psi(\lambda_{n}))^2(c^2+\kappa)},
\end{equation}
\begin{equation}\label{ratefk2}
  R_f^2=\frac{c^2(\kappa-1)^2(c^2+\kappa)\Psi^2(\lambda_{1})\Psi^2(\lambda_{n})}
  {(c^2\Psi(\lambda_{n})+\kappa\Psi(\lambda_{1}))^2(\Psi^2(\lambda_{1})+c^2\kappa\Psi^2(\lambda_{n}))}.
\end{equation}
In addition, if $\Psi(\lambda_{n})=\Psi(\lambda_{1})$ or $c^2=\Psi(\lambda_{1})/\Psi(\lambda_{n})$, then $R_f^1=R_f^2$.
\end{theorem}
\begin{proof}
Let $\epsilon_k=x_k-x^*$. Since $g_k=A\epsilon_k$, by \eqref{gmu}, we have
  \[
    \epsilon_k=\sum_{i=1}^n\lambda_i^{-1}\mu_k^{(i)}\xi_i.
  \]
By Theorem \ref{th1}, we only need to consider the case $\mu_k^{(i)}=0$, $i=2,\ldots,n-1$, that is,
    \[
    \epsilon_k=\lambda_1^{-1}\mu_k^{(1)}\xi_1+\lambda_n^{-1}\mu_k^{(n)}\xi_n.
  \]
  Thus,
    \begin{align}\label{fk}
    f(x_k)-f^*&=\frac{1}{2}\epsilon_k \tr A\epsilon_k
    %=\frac{1}{2}((\mu_k^{(1)})^2\lambda_1^{-1}+(\mu_k^{(n)})^2\lambda_n^{-1})
    %\nonumber\\
    =\frac{1}{2}\frac{\lambda_n(\mu_k^{(1)})^2+\lambda_1(\mu_k^{(n)})^2}{\lambda_1\lambda_n}.
  \end{align}
   Since
      \[
    g_k=\mu_k^{(1)}\xi_1+\mu_k^{(n)}\xi_n \quad \mbox{and}  \quad  \alpha_k=\frac{\Psi(\lambda_{1})(\mu_k^{(1)})^2+\Psi(\lambda_{n})(\mu_k^{(n)})^2}
    {\lambda_1\Psi(\lambda_{1})(\mu_k^{(1)})^2+\lambda_n\Psi(\lambda_{n})(\mu_k^{(n)})^2},
  \]
by the definition of $\epsilon_k$ and the update rule \eqref{eqitr}, we further have that
    \begin{align*}
    \epsilon_{k+1}&=\epsilon_k-\alpha_kg_k
    =(\lambda_1^{-1}-\alpha_k)\mu_k^{(1)}\xi_1+(\lambda_n^{-1}-\alpha_k)\mu_k^{(n)}\xi_n\\
  &=\frac{\Psi(\lambda_{n})(\lambda_n-\lambda_1)(\mu_k^{(n)})^2\mu_k^{(1)}}
  {\lambda_1\left(\lambda_1\Psi(\lambda_{1})(\mu_k^{(1)})^2+\lambda_n\Psi(\lambda_{n})(\mu_k^{(n)})^2\right)}\xi_1\\
  &+\frac{\Psi(\lambda_{1})(\lambda_1-\lambda_n)(\mu_k^{(1)})^2\mu_k^{(n)}}
  {\lambda_n\left(\lambda_1\Psi(\lambda_{1})(\mu_k^{(1)})^2+\lambda_n\Psi(\lambda_{n})(\mu_k^{(n)})^2\right)}\xi_n\\
  &=
  \frac{(\lambda_n-\lambda_1)\left(\lambda_n\Psi(\lambda_{n})(\mu_k^{(n)})^2\mu_k^{(1)}\xi_1
  -\lambda_1\Psi(\lambda_{1})(\mu_k^{(1)})^2\mu_k^{(n)}\xi_n\right)}
  {\lambda_1\lambda_n\left(\lambda_1\Psi(\lambda_{1})(\mu_k^{(1)})^2+\lambda_n\Psi(\lambda_{n})(\mu_k^{(n)})^2\right)}.
  %&=\frac{(\lambda_n-\lambda_1)
%  (\lambda_n^{u+1}\mu_k^{(1)}(\mu_k^{(n)})^2\xi_1-\lambda_1^{u+1}(\mu_k^{(1)})^2\mu_k^{(n)}\xi_n)}
%  {\lambda_1\lambda_n(\lambda_1^{u+1}(\mu_k^{(1)})^2+\lambda_n^{u+1}(\mu_k^{(n)})^2)}.
  \end{align*}
Hence, we obtain
\begin{align}\label{fk1}
  &f(x_{k+1})-f^*
  =\frac{1}{2}\epsilon_{k+1} \tr A\epsilon_{k+1}\nonumber\\
  =&
  \frac{1}{2}\frac{(\lambda_n-\lambda_1)^2(\mu_k^{(1)})^2(\mu_k^{(n)})^2
  \left(\lambda_n\Psi^2(\lambda_{n})(\mu_k^{(n)})^2+\lambda_1\Psi^2(\lambda_{1})(\mu_k^{(1)})^2\right)}
  {\lambda_1\lambda_n\left(\lambda_1\Psi(\lambda_{1})(\mu_k^{(1)})^2+\lambda_n\Psi(\lambda_{n})(\mu_k^{(n)})^2\right)^2}.
\end{align}
Combining \eqref{fk} with \eqref{fk1} yields that
\begin{align*}
&  \frac{f(x_{k+1})-f^*}{f(x_{k})-f^*} =\frac{\epsilon_{k+1} \tr A\epsilon_{k+1}}{\epsilon_k \tr A\epsilon_k}\\
%  &=
%\frac{(\lambda_n-\lambda_1)^2(\mu_k^{(1)})^2(\mu_k^{(n)})^2
%  \left(\lambda_n\Psi^2(\lambda_{n})(\mu_k^{(n)})^2+\lambda_1\Psi^2(\lambda_{1})(\mu_k^{(1)})^2\right)}
%  {\left(\lambda_1\Psi(\lambda_{1})(\mu_k^{(1)})^2+\lambda_n\Psi(\lambda_{n})(\mu_k^{(n)})^2\right)^2
%  \left(\lambda_n(\mu_k^{(1)})^2+\lambda_1(\mu_k^{(n)})^2\right)}\\
%  &=\frac{(\lambda_n-\lambda_1)^2(\mu_k^{(1)})^2(\mu_k^{(n)})^2
%  (\lambda_n^2\Psi^2(\lambda_{n})(\mu_k^{(n)})^2+\lambda_1^2\Psi^2(\lambda_{1})(\mu_k^{(1)})^2)}
%  {\lambda_1\lambda_n(\lambda_1\Psi(\lambda_{1})(\mu_k^{(1)})^2+\lambda_n\Psi(\lambda_{n})(\mu_k^{(n)})^2)^2}
%  \frac{1}{(\lambda_n(\mu_k^{(1)})^2+\lambda_1(\mu_k^{(n)})^2)}\\
  = &\frac{(\mu_k^{(1)})^2(\mu_k^{(n)})^2(\kappa-1)^2
  \left(\kappa\Psi^2(\lambda_{n})(\mu_k^{(n)})^2+\Psi^2(\lambda_{1})(\mu_k^{(1)})^2\right)}
  {\left(\Psi(\lambda_{1})(\mu_k^{(1)})^2+\kappa\Psi(\lambda_{n})(\mu_k^{(n)})^2\right)^2
  \left(\kappa(\mu_k^{(1)})^2+(\mu_k^{(n)})^2\right)},
\end{align*}
  which gives \eqref{f2k2f1k} by substituting the limits of $(\mu_k^{(1)})^2$ and $(\mu_k^{(n)})^2$ in Theorem \ref{th1}.

Notice $\kappa>1$ by our assumption. So,  $R_f^1=R_f^2$ is equivalent to
  \[
  \frac{\Psi^2(\lambda_{1})+c^2\kappa\Psi^2(\lambda_{n})}
  {(\Psi(\lambda_{1})+c^2\kappa\Psi(\lambda_{n}))^2(c^2+\kappa)}
=\frac{(c^2+\kappa)\Psi^2(\lambda_{1})\Psi^2(\lambda_{n})}
  {(c^2\Psi(\lambda_{n})+\kappa\Psi(\lambda_{1}))^2(\Psi^2(\lambda_{1})+c^2\kappa\Psi^2(\lambda_{n}))},
\]
%i.e.,
%  \begin{eqnarray*}
%   && (c^2\Psi(\lambda_{n})+\kappa\Psi(\lambda_{1}))^2(\Psi^2(\lambda_{1})+c^2\kappa\Psi^2(\lambda_{n}))^2 \\
%  &=& (c^2+\kappa)^2\Psi^2(\lambda_{1})\Psi^2(\lambda_{n})(\Psi(\lambda_{1})+c^2\kappa\Psi(\lambda_{n}))^2.
%\end{eqnarray*}
which by rearranging terms gives
  \[
  c^4\Psi^2(\lambda_{n})(\Psi(\lambda_{n})-\Psi(\lambda_{1}))=\Psi^2(\lambda_{1})(\Psi(\lambda_{n})-\Psi(\lambda_{1})).
\]
Hence,   $R_f^1=R_f^2$  holds if $\Psi(\lambda_{n})=\Psi(\lambda_{1})$ or $c^2=\Psi(\lambda_{1})/\Psi(\lambda_{n})$.
\end{proof}

\begin{remark}
Theorem \ref{thconfk} indicates that, when $\Psi(A)=I$ (i.e., the SD method), the two sequences
$\Big\{\frac{\Delta_{2k+1}}{\Delta_{2k}} \Big\}$ and $\Big\{ \frac{\Delta_{2k+2}}{\Delta_{2k+1}} \Big\}$
converge at the same speed, where $\Delta_k = f(x_k) - f^*$.
Otherwise, the two sequences may converge at different rates.
\end{remark}

\begin{figure}[h]
  \centering
  \includegraphics[width=0.75\textwidth,height=0.52\textwidth]{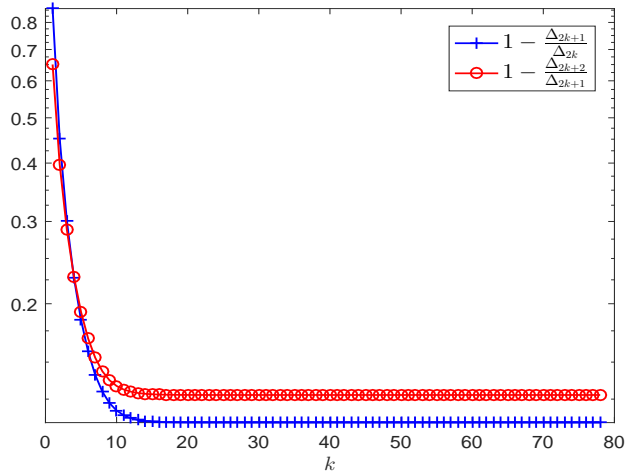}\\
  \caption{Problem \eqref{exfkrate} with $n=10$: convergence history of the sequences $\big\{1- \frac{\Delta_{2k+1}}{\Delta_{2k}}\big\}$
and $\big\{1- \frac{\Delta_{2k+2}}{\Delta_{2k+1}}\big\}$
  generated by gradient method \eqref{glstep1} with $\Psi(A)=A$  (i.e., the MG method).}\label{confk}
\end{figure}

To illustrate the results in Theorem \ref{thconfk}, we apply gradient method \eqref{glstep1} with $\Psi(A)=A$ (i.e., the MG method) to an instance of \eqref{qudpro},
where the vector of all ones was used as the initial point, the matrix $A$ is diagonal with
\begin{equation}\label{exfkrate}
  A_{ii}=i\sqrt{i},~~i=1,\ldots,n,
\end{equation}
 and $b=0$. Figure \ref{confk} clearly shows the difference between $R_f^1$ and $R_f^2$.

The next theorem shows the asymptotic convergence of the gradient norm.
\begin{theorem}\label{thcongk}
Under the conditions of Theorem \ref{th1}, the following limits hold,
  \begin{equation}\label{gk2g1}
    \lim_{k\rightarrow\infty}\frac{\|g_{2k+1}\|^2}{\|g_{2k}\|^2}
     =R_g^1
\quad \mbox{ and } \quad     \lim_{k\rightarrow\infty}\frac{\|g_{2k+2}\|^2}{\|g_{2k+1}\|^2}
     =R_g^2,
  \end{equation}
where
\begin{equation}\label{rg}
  R_g^1=\frac{c^2(\kappa-1)^2(\Psi^2(\lambda_{1})+c^2\Psi^2(\lambda_{n}))}{(1+c^2)(\Psi(\lambda_{1})+c^2\kappa\Psi(\lambda_{n}))^2},
\end{equation}
\begin{equation}\label{rg2}
  R_g^2=\frac{c^2(1+c^2)(\kappa-1)^2\Psi^2(\lambda_{1})\Psi^2(\lambda_{n})}
  {(c^2\Psi(\lambda_{n})+\kappa\Psi(\lambda_{1}))^2(\Psi^2(\lambda_{1})+c^2\Psi^2(\lambda_{n}))}.
\end{equation}
In addition, if $\Psi(\lambda_{n})=\kappa\Psi(\lambda_{1})$ or $c^2=\Psi(\lambda_{1})/\Psi(\lambda_{n})$, then $R_g^1=R_g^2$.
\end{theorem}
\begin{proof}
Using the same arguments as in Theorem \ref{thconfk}, we have
\[
    \|g_{k}\|^2=(\mu_k^{(1)})^2+(\mu_k^{(n)})^2
  \]
and
%\begin{align*}
\[
  \|g_{k+1}\|^2=\epsilon_{k+1} \tr A^2\epsilon_{k+1}
  =
  \frac{(\lambda_n-\lambda_1)^2(\mu_k^{(1)})^2(\mu_k^{(n)})^2
  \left(\Psi^2(\lambda_{n})(\mu_k^{(n)})^2+\Psi^2(\lambda_{1})(\mu_k^{(1)})^2\right)}
  {\left(\lambda_1\Psi(\lambda_{1})(\mu_k^{(1)})^2+\lambda_n\Psi(\lambda_{n})(\mu_k^{(n)})^2\right)^2},
\]
%\end{align*}
which give that
\[
  \frac{\|g_{k+1}\|^2}{\|g_{k}\|^2}
%  &\frac{(\lambda_n-\lambda_1)^2(\mu_k^{(1)})^2(\mu_k^{(n)})^2
%  (\lambda_n^{2u}(\mu_k^{(n)})^2+\lambda_1^{2u}(\mu_k^{(1)})^2)}
%  {(\lambda_1^{u+1}(\mu_k^{(1)})^2+\lambda_n^{u+1}(\mu_k^{(n)})^2)^2((\mu_k^{(1)})^2+(\mu_k^{(n)})^2)}\\
 =
    \frac{(\kappa-1)^2(\mu_k^{(1)})^2(\mu_k^{(n)})^2
  \left(\Psi^2(\lambda_{n})(\mu_k^{(n)})^2+\Psi^2(\lambda_{1})(\mu_k^{(1)})^2\right)}
  {\left(\Psi(\lambda_{1})(\mu_k^{(1)})^2+\kappa\Psi(\lambda_{n})(\mu_k^{(n)})^2\right)^2
  \left((\mu_k^{(1)})^2+(\mu_k^{(n)})^2\right)}.
\]
Thus, \eqref{gk2g1} follows by substituting the limits of $(\mu_k^{(1)})^2$ and $(\mu_k^{(n)})^2$ in Theorem \ref{th1}.

Notice $\kappa>1$ by our assumption. So,  $R_g^1=R_g^2$ is equivalent to
  \[
 \frac{\Psi^2(\lambda_{1})+c^2\Psi^2(\lambda_{n})}{(1+c^2)(\Psi(\lambda_{1})+c^2\kappa\Psi(\lambda_{n}))^2}
=\frac{(1+c^2)\Psi^2(\lambda_{1})\Psi^2(\lambda_{n})}
  {(c^2\Psi(\lambda_{n})+\kappa\Psi(\lambda_{1}))^2
  (\Psi^2(\lambda_{1})+c^2\Psi^2(\lambda_{n}))},
\]
%i.e.,
%  \begin{eqnarray*}
% && (c^2\Psi(\lambda_{n})+\kappa\Psi(\lambda_{1}))^2(\Psi^2(\lambda_{1})+c^2\Psi^2(\lambda_{n}))^2  \\
%  &= &
%  (1+c^2)^2\Psi^2(\lambda_{1})\Psi^2(\lambda_{n})(\Psi(\lambda_{1})+c^2\kappa\Psi(\lambda_{n}))^2.
%\end{eqnarray*}
which by rearranging terms gives
  \[
  c^4\Psi^2(\lambda_{n})(\kappa\Psi(\lambda_{1})-\Psi(\lambda_{n}))=
  \Psi^2(\lambda_{1})(\kappa\Psi(\lambda_{1})-\Psi(\lambda_{n})).
\]
Hence,  $R_g^1=R_g^2$ holds if $\Psi(\lambda_{n})=\kappa\Psi(\lambda_{1})$ or $c^2=\Psi(\lambda_{1})/\Psi(\lambda_{n})$.
\end{proof}

\begin{remark}\label{gratemk}
Theorem \ref{thcongk}  indicates that the two sequences $\Big\{\frac{\|g_{2k+1}\|^2}{\|g_{2k}\|^2}\Big\}$ and \newline $\Big\{\frac{\|g_{2k+2}\|^2}{\|g_{2k+1}\|^2}\Big\}$ generated by the MG method (i.e., $\Psi(A)=A$) converge at the same rate. Otherwise, the two sequences may converge at different rates.
\end{remark}

By Theorems \ref{thconfk} and \ref{thcongk}, we can obtain the following corollary.
%showing that the odd and even subsequences of objective values/gradient norms converge at the same rate.
\begin{corollary}\label{ratefg}
Under the conditions of Theorem \ref{th1}, we have
    \begin{equation}\label{ratefeq}
     \lim_{k\rightarrow\infty}\frac{f(x_{2k+3})-f^*}{f(x_{2k+1})-f^*}
     =\lim_{k\rightarrow\infty}\frac{f(x_{2k+2})-f^*}{f(x_{2k})-f^*}
     =R_f^1R_f^2,
     %=\frac{c^4(\kappa-1)^4}{(1+c^2\kappa^{u+1})^2(c^2+\kappa^{1-u})^2}.
     %=R_g^1R_g^2
  \end{equation}
      \begin{equation}\label{rategeq}
    \lim_{k\rightarrow\infty}\frac{\|g_{2k+3}\|^2}{\|g_{2k+1}\|^2}
     =\lim_{k\rightarrow\infty}\frac{\|g_{2k+2}\|^2}{\|g_{2k}\|^2}
     =R_g^1R_g^2.
     %=\frac{c^4(\kappa-1)^4}{(1+c^2\kappa^{u+1})^2(c^2+\kappa^{1-u})^2}.
     %=R_g^1R_g^2
  \end{equation}
  In addition,
  \begin{equation}\label{ratefgeq}
    R_f^1R_f^2=R_g^1R_g^2=\frac{c^4(\kappa-1)^4\Psi^2(\lambda_{1})\Psi^2(\lambda_{n})}
    {(\Psi(\lambda_{1})+c^2\kappa\Psi(\lambda_{n}))^2(c^2\Psi(\lambda_{n})+\kappa\Psi(\lambda_{1}))^2}.
  \end{equation}
\end{corollary}

\begin{remark}\label{fgrate}
Corollary \ref{ratefg} shows that the odd and even subsequences of objective values and gradient norms converge at the same rate. Moreover, we have
  \begin{equation}\label{ratefgupbd}
    R_f^1R_f^2=R_g^1R_g^2
    =\frac{(\kappa-1)^4}
    {(1+\kappa/t+t\kappa+\kappa^2)^2}
    \leq\left(\frac{\kappa-1}{\kappa+1}\right)^4,
  \end{equation}
where $t= c^2 \Psi(\lambda_{n})/\Psi(\lambda_{1})$.
Notice that the right side of \eqref{ratefgupbd} only depends on $\kappa$, which implies these
 odd and even subsequences  generated by all the gradient methods \eqref{glstep1} will have the same worst asymptotic rate
independent of $\Psi$.
\end{remark}

%We close this section by deriving a bound on the constant $c$ defined in Theorem \ref{th1}. As in \cite{nocedal2002behavior}, we employ the \emph{minimum deviation}
Now, as in \cite{nocedal2002behavior}, we define the \emph{minimum deviation}
\begin{equation}\label{defmsig}
  \sigma=\min_{i\in\mathcal{I}}~\left| \frac{2 \lambda_i- (\lambda_1+\lambda_n)}{ \lambda_n-\lambda_1 }\right|,
\end{equation}
where
\begin{equation*}
  \mathcal{I}=\{i: \lambda_1<\lambda_i<\lambda_n,~g_0 \tr \xi_i\neq0,~\textrm{and}~\lambda_i\neq\alpha_k~\textrm{for all}~k\}.
\end{equation*}
Clearly, $\sigma\in(0,1)$. We now close this section by deriving a bound on the constant $c$ defined in Theorem \ref{th1}.
The following theorem generalizes the results in \cite{akaike1959successive,nocedal2002behavior}, where only the case $\Psi(A)=I$ (i.e., the SD method) is considered.
\begin{theorem}\label{bdc2}
 Under the conditions of Theorem \ref{th1}, and assuming that $\mathcal{I}$ is nonempty, we have
\begin{equation}\label{eqbdc2}
  \frac{\Psi(\lambda_{1})}{\Psi(\lambda_{n})}\frac{1}{\phi_\sigma}\leq c^2\leq\frac{\Psi(\lambda_{1})}{\Psi(\lambda_{n})}\phi_\sigma,
\end{equation}
where
\begin{equation}\label{defphi}
  \phi_\sigma=\frac{2+\eta_\sigma+\sqrt{\eta_\sigma^2+4\eta_\sigma}}{2} \quad \mbox{and} \quad \eta_\sigma=4\left(\frac{1+\sigma^2}{1-\sigma^2}\right).
\end{equation}
\end{theorem}
\begin{proof}
%  Suppose that $\lim_{k\rightarrow\infty}T^{k}p=p_*$, by Lemma \ref{tpeqp} and Lemma \ref{lm3}, we have that $\lim_{k\rightarrow\infty}T^{k+1}p=Tp_*$ and $\lim_{k\rightarrow\infty}T^{k+2}p=p_*$.with $q_0$ defined by \eqref{defpk}
  Let $p=q_0$. By the definition of $T$, we have that
      \begin{equation}\label{tpk2i21}
        \frac{(T^{k+2}p)^{(i)}}{(T^{k+2}p)^{(1)}}=
    \frac{(T^{k}p)^{(i)}}{(T^{k}p)^{(1)}}
    \frac{(\lambda_i-\gamma(T^{k}p))^2(\lambda_i-\gamma(T^{k+1}p))^2}
  {(\lambda_1-\gamma(T^{k}p))^2(\lambda_1-\gamma(T^{k+1}p))^2}.
  \end{equation}
It follows from Theorem \ref{th1} and Lemma \ref{lm3} that
  \begin{equation}\label{tpki}
    \frac{(T^{k}p)^{(i)}}{(T^{k}p)^{(1)}}\rightarrow0,~i=2,\ldots,n-1.
  \end{equation}
   By the continuity of $T$ and \eqref{t2k} in Lemma \ref{lm3}, we always have that
      \begin{equation*}\label{eq67}
  \frac{(\lambda_i-\gamma(T^{k}p))^2(\lambda_i-\gamma(T^{k+1}p))^2}
  {(\lambda_1-\gamma(T^{k}p))^2(\lambda_1-\gamma(T^{k+1}p))^2}\rightarrow
    \frac{(\lambda_i-\gamma(p_*))^2(\lambda_i-\gamma(Tp_*))^2}
  {(\lambda_1-\gamma(p_*))^2(\lambda_1-\gamma(Tp_*))^2},
  \end{equation*}
which together with \eqref{tpk2i21} and \eqref{tpki} implies that
% \begin{equation}\label{tpk21}, for sufficiently large $k$,
%   \frac{(\lambda_i-\gamma(T^{k}p))^2(\lambda_i-\gamma(T^{k+1}p))^2}
%  {(\lambda_1-\gamma(T^{k}p))^2(\lambda_1-\gamma(T^{k+1}p))^2}\leq1,~i=2,\ldots,n-1.
% \end{equation}
%Combining \eqref{eq67} with \eqref{tpk21}, we obtain
    \begin{equation}\label{lbdgama1}
  \frac{(\lambda_i-\gamma(p_*))^2(\lambda_i-\gamma(Tp_*))^2}
  {(\lambda_1-\gamma(p_*))^2(\lambda_1-\gamma(Tp_*))^2}
  \leq1,~i=2,\ldots,n-1,
  \end{equation}
where $p_*$ is the same vector as in Lemma \ref{lm3}. Clearly, \eqref{lbdgama1} also holds for $i=1$.
As for $i=n$, it follows from \eqref{gmaptp} in Lemma \ref{tpeqp} and Theorem \ref{th1} that
\begin{equation}\label{gmaptp7}
  \gamma(p_*)+\gamma(Tp_*)=\lambda_1+\lambda_n,
\end{equation}
  which yields that
  \[
  \frac{(\lambda_n-\gamma(p_*))^2(\lambda_n-\gamma(Tp_*))^2}
  {(\lambda_1-\gamma(p_*))^2(\lambda_1-\gamma(Tp_*))^2}=1.
  \]
Thus, \eqref{lbdgama1} holds for $i=1,\ldots,n$. Hence, we have
 \begin{eqnarray}\label{ineqsig2}
 && \left(\lambda_i-\delta
  -\left(\gamma(p_*)-\delta\right)\right)^2
  \left(\lambda_i-\delta
  -\left(\gamma(Tp_*)-\delta\right)\right)^2
  \\
& \leq & \left(\lambda_1-\delta
-\left(\gamma(p_*)-\delta\right)\right)^2
  \left(\lambda_1-\delta
  -\left(\gamma(Tp_*)-\delta\right)\right)^2, \nonumber
  \end{eqnarray}
 where $\delta=\frac{\lambda_1+\lambda_n}{2}$.
  By \eqref{gmaptp7} and \eqref{ineqsig2}, we obtain
\begin{eqnarray*}
&&  \left(\lambda_i-\delta
  -\left(\gamma(p_*)-\delta\right)\right)^2
  \left(\lambda_i-\delta
  +\left(\gamma(p_*)-\delta\right)\right)^2 \\
&\leq & \left(\frac{\lambda_1-\lambda_n}{2}
-\left(\gamma(p_*)-\delta\right)\right)^2
  \left(\frac{\lambda_1-\lambda_n}{2}
  +\left(\gamma(p_*)-\delta\right)\right)^2,
  \end{eqnarray*}
%  i.e.,
%  \begin{align*}
%  &\left(\left(\lambda_i-\delta\right)^2
%  -\left(\gamma(p_*)-\delta\right)^2\right)^2
% \leq\left(\left(\frac{\lambda_1-\lambda_n}{2}\right)^2
% -\left(\gamma(p_*)-\delta\right)^2\right)^2,
%  \end{align*}
  which implies that
  \begin{equation}\label{ineqsig1}
    \left(\frac{\lambda_1-\lambda_n}{2}\right)^2+
    \left(\lambda_i-\delta\right)^2
    \geq
    2\left(\gamma(p_*)-\delta\right)^2.
  \end{equation}
By Lemma \ref{tpeqp} and Theorem \ref{th1}, we have that
    \[
    \gamma(p_*)=\frac{\lambda_1\Psi(\lambda_{1})p_*^{(1)}+\lambda_n\Psi(\lambda_{n})p_*^{(n)}}
    {\Psi(\lambda_{1})p_*^{(1)}+\Psi(\lambda_{n})p_*^{(n)}}.
  \]
Substituting $\gamma(p_*)$ into \eqref{ineqsig1}, we obtain
    \begin{equation*}\label{ineqsig3}
    \left(\frac{\lambda_1-\lambda_n}{2}\right)^2+
    \left(\lambda_i-\delta\right)^2
    \geq
     \frac{(\lambda_n-\lambda_1)^2(\Psi(\lambda_{n})c^2-\Psi(\lambda_{1}))^2}
  {2(\Psi(\lambda_{n})c^2+\Psi(\lambda_{1}))^2},
  \end{equation*}
  which gives
  \begin{equation}\label{inqsig}
    4\left(\frac{1+\sigma_i^2}{1-\sigma_i^2}\right)\geq\frac{(c^2\Psi(\lambda_{n})-\Psi(\lambda_{1}))^2}{c^2\Psi(\lambda_{1})\Psi(\lambda_{n})},
\quad \mbox{where} \quad   \sigma_i=\frac{2 \lambda_i-(\lambda_1+\lambda_n)}
  {\lambda_n-\lambda_1}.
  \end{equation}
Noting that \eqref{inqsig} holds for all $i\in\mathcal{I}$. Thus, we have
  \begin{equation}\label{inqsig2}
    \frac{(c^2\Psi(\lambda_{n})-\Psi(\lambda_{1}))^2}{c^2\Psi(\lambda_{1})\Psi(\lambda_{n})}\leq\eta_\sigma,
  \end{equation}
  which implies \eqref{eqbdc2}. This completes the proof.
\end{proof}

\section{Techniques for breaking the zigzagging pattern}\label{secnewstep}

As shown in the previous section, all the gradient methods \eqref{glstep1} asymptotically conduct its searches
 in the two-dimensional subspace spanned by $\xi_1$ and $\xi_n$. By \eqref{upmuk}, if either
 $\mu_k^{(1)}$ or $\mu_k^{(n)}$ equals to zero, the corresponding component will vanish at all subsequent iterations. Hence, in order to break the undesired zigzagging pattern, a good strategy is to
employ some stepsize approximating $1/\lambda_1$ or $1/\lambda_n$. In this section, we will derive a
new stepsize converging to $1/\lambda_n$ and propose a periodic gradient method using this new
stepsize.

\subsection{A new stepsize}\label{subsecns}
Our new stepsize will be derived by imposing finite termination on minimizing two-dimensional strictly convex quadratic function, see \cite{yuan2006new} for the case of $\Psi(A)=I$ (i.e., the SD method).
We mention that the key property used by Yuan \cite{yuan2006new} is that two consecutive gradients
generated by the SD method are perpendicular to each other, which may not be true for all the gradient methods
\eqref{glstep1}. However, we have by the stepsize definition \eqref{glstep1}  that
\begin{equation}\label{orthg}
  g_k \tr \Psi(A)g_{k+1}=g_k \tr \Psi(A)g_k-\alpha_kg_k \tr \Psi(A)Ag_k=0.
\end{equation}

Consider the two-dimensional case. Suppose we want to find the
minimizer of \eqref{qudpro} with $n=2$ after the following $3$ iterations:
\[
  x_1 = x_0 - \alpha_0g_0, \quad
  x_2 = x_1 - \alpha_1g_1, \quad
  x_3 = x_2 - \alpha_2g_2,
\]
where $g_i \ne 0$, $i=0,1,2$,
 $\alpha_0$ and $\alpha_2$ are stepsizes given by \eqref{glstep1}, and $\alpha_1$ is to
be derived by ensuring $x_3$ is the solution.

By \eqref{orthg}, we have $g_0 \tr \Psi(A)g_1=0$. Hence, all vectors $x_k$ can be expressed by the linear combination of $\frac{\Psi^r(A)g_0}{\|\Psi^r(A)g_0\|}$ and $\frac{\Psi^{1-r}(A)g_1}{\|\Psi^{1-r}(A)g_1\|}$ for any given $r \in \mathbb{R}$. Now, consider
\begin{eqnarray}\label{ftu1}
\varphi(t,l) &:= &
 f\left(x_1+t\frac{\Psi^r(A)g_0}{\|\Psi^r(A)g_0\|}+l\frac{\Psi^{1-r}(A)g_1}{\|\Psi^{1-r}(A)g_1\|}\right) \\
& = &
f(x_1) +G \tr
              \begin{pmatrix}
                t \\
                l \\
              \end{pmatrix}
            +  \frac{1}{2}
              \begin{pmatrix}
                t \\
                l \\
              \end{pmatrix} \tr
H
              \begin{pmatrix}
                t \\
                l \\
              \end{pmatrix}, \nonumber
\end{eqnarray}
where
\begin{equation}\label{eqgrad}
  G= B g_1 =  \begin{pmatrix}
                \frac{g_1 \tr \Psi^r(A)g_0}{\|\Psi^r(A)g_0\|}\\
                \frac{g_1 \tr \Psi^{1-r}(A)g_1}{\|\Psi^{1-r}(A)g_1\|}
              \end{pmatrix}
\mbox{ with  }    B=\begin{pmatrix}
      \frac{\Psi^r(A)g_0}{\|\Psi^r(A)g_0\|},
      \frac{\Psi^{1-r}(A)g_1}{\|\Psi^{1-r}(A)g_1\|} \\
    \end{pmatrix} \tr
\end{equation}
and
\begin{equation} \label{eqhes}
  H=  B A B \tr  =  \begin{pmatrix}
      \frac{g_0 \tr \Psi^{2r}(A)Ag_0}{\|\Psi^r(A)g_0\|^2} & \frac{g_0 \tr \Psi(A) Ag_1}{\|\Psi^r(A)g_0\|\|\Psi^{1-r}(A)g_1\|} \\
      \frac{g_0 \tr \Psi(A)Ag_1}{\|\Psi^r(A)g_0\|\|\Psi^{1-r}(A)g_1\|}  & \frac{g_1 \tr \Psi^{2(1-r)}(A)Ag_1}{\|\Psi^{1-r}(A)g_1\|^2}\\
    \end{pmatrix}.
\end{equation}
Note that $B \tr B = B B \tr = I$ since $n=2$.
The minimizer $(t^*, l^*)$ of $\varphi$ satisfy
\begin{equation*}
G
            +
    H
              \begin{pmatrix}
                t^* \\
                l^* \\
              \end{pmatrix}=0,
\quad \Longrightarrow \quad
              \begin{pmatrix}
                t^* \\
                l^* \\
              \end{pmatrix}
=-H^{-1}G.
%= -\frac{1}{\Delta}\hat{H}G,
\end{equation*}
Suppose $x_3$ is the solution, that is
\begin{equation*}
  x_3=x_1+t^*\frac{\Psi^r(A)g_0}{\|\Psi^r(A)g_0\|}+l^*\frac{\Psi^{1-r}(A)g_1}{\|\Psi^{1-r}(A)g_1\|}.
\end{equation*}
Then, since $x_3 = x_2 - \alpha_2 g_2$, we have  $x_3-x_2$ is parallel to $g_2$, i.e.,
\begin{equation} \label{parallel-g2}
B \tr  \begin{pmatrix}
                t^* \\
                l^* \\
              \end{pmatrix}
+ \alpha_1 g_1
\quad \mbox{ is parallel to } \quad g_2,
\end{equation}
which is equivalent to
\begin{equation}\label{g2pal1}
              \begin{pmatrix}
                t^* \\
                l^* \\
              \end{pmatrix}
            -(-\alpha_1G)
            = -(H^{-1}G-\alpha_1 G)
\quad \mbox{and} \quad  G+H(-\alpha_1G)
\end{equation}
are parallel. Denote the components of $G$ by $G_i$, and the components of $H$ by $H_{ij}$, $i,j=1,2$.
By \eqref{g2pal1}, we would have
\begin{equation*}
              \begin{pmatrix}
                H_{22}G_1-H_{12}G_2-\alpha_1\Delta G_1 \\
                H_{11}G_2-H_{12}G_1-\alpha_1\Delta G_2\\
              \end{pmatrix}
 \quad \mbox{and} \quad
              \begin{pmatrix}
                G_1-\alpha_1(H_{11}G_1+H_{12}G_2) \\
                G_2-\alpha_1(H_{12}G_1+H_{22}G_2)\\
              \end{pmatrix}.
\end{equation*}
are parallel, where $\Delta=\textrm{det}(H) = \textrm{det}(A)  >0$. It follows that
\begin{eqnarray*}
 && (H_{22}G_1-H_{12}G_2-\alpha_1\Delta G_1)
  [G_2-\alpha_1(H_{12}G_1+H_{22}G_2)]\\
 & =& (H_{11}G_2-H_{12}G_1-\alpha_1\Delta G_2)
  [G_1-\alpha_1(H_{11}G_1+H_{12}G_2)],
\end{eqnarray*}
which gives
\begin{equation}\label{eqalp1}
  \alpha_1^2\Delta \Gamma-\alpha_1(H_{11}+H_{22})\Gamma+\Gamma=0,
\end{equation}
where
\[
  \Gamma=(H_{12}G_1+H_{22}G_2)G_1-(H_{11}G_1+H_{12}G_2)G_2.
\]
On the other hand, if \eqref{eqalp1} holds, we have \eqref{parallel-g2} holds,
which by \eqref{eqgrad},  $H^{-1} = B A^{-1} B \tr$ and $B \tr B=I$ implies that
\[
- B \tr  H^{-1} G + \alpha_1 g_1 = -  A^{-1} g_1 + \alpha_1 g_1 = - A^{-1}(g_1- \alpha_1 A g_1)
= - A^{-1} g_2
\]
is parallel to $g_2$. Hence, $g_2$ is an eigenvector of $A$, i.e. $Ag_2 = \lambda g_2$ for some
$\lambda > 0$, since $g_2 \ne 0$.
So, by \eqref{glstep1}, $\alpha_2 = \Psi(\lambda) g_2 \tr  g_2 /(\lambda \Psi(\lambda) g_2 \tr  g_2  ) = 1/\lambda$.
Therefore, $g_3 = g_2 - \alpha_2 A g_2 = g_2 - \alpha_2 \lambda g_2 = 0$, which implies $x_3$ is
the solution. So, \eqref{eqalp1} guarantees $x_3$ is the minimizer.

%When $\Gamma=0$, we have $HG$ is parallel to $G$. Hence, since $G \ne 0$, $G$ is an eigenvector of $H$,
%that is, $HG=\lambda G$ for some positive $\lambda$.
%  Then, we have by \eqref{eqgrad} that
%\begin{equation*}
% HB g_1= HG=\lambda G=\lambda B g_1,
%\end{equation*}
%which yields $B \tr HBg_1=\lambda g_1$, i.e., $Ag_1=\lambda g_1$.
%Hence, in this case $g_1$ is an eigenvector of $A$.
%So, if $\alpha_1=1/\lambda$, $x_2$ is the solution; otherwise, $x_3$ must be the solution.

Hence, to ensure $x_3$ is the minimizer, by \eqref{eqalp1}, we only need
to choose $\alpha_1$ satisfying
\begin{equation}\label{eqalp2}
  \alpha_1^2\Delta -\alpha_1(H_{11}+H_{22})+1=0,
\end{equation}
whose two positive roots are
\begin{equation*}
  \frac{(H_{11}+H_{22})\pm\sqrt{(H_{11}+H_{22})^2-4\Delta}}{2\Delta}.
\end{equation*}
These two roots are $1/\lambda_1$ and $1/\lambda_2$, where $0 < \lambda_1 < \lambda_2$
are two eigenvalues of $A$ (Note that $A$ and $H$ have same eigenvalues).
For numerical reasons (see next subsection), we would like to choose $\alpha_1$ to
be the smaller one $1/\lambda_2$, which can be calculated as
\begin{align}\label{new2}
  \alpha_1&= \frac{2}{(H_{11}+H_{22})+\sqrt{(H_{11}+H_{22})^2-4\Delta}}
  \nonumber\\
  &= \frac{2}{(H_{11}+H_{22})+\sqrt{(H_{11}-H_{22})^2+4H_{12}^2}}.
\end{align}

To check this finite termination property, we applied the above described method with $\alpha_1$ given by \eqref{new2}, and $\Psi(A)=A$ in \eqref{glstep1}, (i.e., $\alpha_0$ and $\alpha_2$
use the MG stepsize) to minimize two-dimensional quadratic function \eqref{qudpro} with
\begin{equation}\label{twoquad}
  A=\textrm{diag}\{1,\lambda\} \quad \mbox{and} \quad b=0.
\end{equation}
%We stopped the iteration if $\|g_k\|\leq10^{-8}$.
We run the algorithm for 3 iterations using ten random starting points and
the averaged values of $\|g_3\|$ and $f(x_3)$ are presented in Table \ref{tb2ft}.
We can observe that for different values of $\lambda$, the $\|g_3\|$ and $f(x_3)$ obtained by
the method in three iterations are numerically very close to zero. This coincides with our analysis.

\begin{table}[ht!b]
\setlength{\tabcolsep}{4ex}
{\footnotesize
\caption{Averaged results for problem \eqref{twoquad} with different condition numbers.}\label{tb2ft}
\begin{center}
\begin{tabular}{|c|c|c|}
\hline
 \multicolumn{1}{|c|}{\multirow{1}{*}{$\lambda$}}
 &\multicolumn{1}{c|}{$\|g_3\|$} &\multicolumn{1}{c|}{$f(x_3)$}\\
 \hline
 10       &   4.8789e-18  &   8.0933e-36\\
 \hline
100       &   4.1994e-18  &   2.2854e-37\\
\hline
1000     &   1.2001e-18  &   2.8083e-39\\
\hline
10000    &   1.0621e-18  &   5.3460e-40\\
\hline
%100000  &3 &   8.4802e-21  &   2.1245e-45\\
%\hline
\end{tabular}
\end{center}
}
\end{table}

\subsection{Spectral property of the new stepsize}\label{astilalp}

Notice that $g_1=g_0-\alpha_0Ag_0$ and $g_0 \tr \Psi(A)g_1=0$. So, we have
\[
 g_0 \tr \Psi(A)Ag_1=-(g_1 \tr \Psi(A)g_1)/\alpha_0.
\]
Hence, the matrix $H$ given in \eqref{eqhes} can be also written as
\begin{equation}\label{eqhes2}
  H=    \begin{pmatrix}
      \frac{g_0 \tr \Psi^{2r}(A)Ag_0}{\|\Psi^r(A)g_0\|^2} & -\frac{g_1 \tr \Psi(A) g_1}{\alpha_0\|\Psi^r(A)g_0\|\|\Psi^{1-r}(A)g_1\|} \\
      -\frac{g_1 \tr \Psi(A) g_1}{\alpha_0\|\Psi^r(A)g_0\|\|\Psi^{1-r}(A)g_1\|}  & \frac{g_1 \tr \Psi^{2(1-r)}(A)Ag_1}{\|\Psi^{1-r}(A)g_1\|^2} \\
    \end{pmatrix}.
\end{equation}
So, for general case,  we could propose our new stepsize at the $k$-th iteration as
\begin{equation}\label{newn}
   \tilde{\alpha}_{k}
  =\frac{2}{(H_{11}^k+H_{22}^k)+\sqrt{(H_{11}^k-H_{22}^k)^2+4(H_{12}^k)^2}},
\end{equation}
where $H_{ij}^k$ is the component of $H^k$:
\begin{equation}\label{eqhesk}
  H^k=    \begin{pmatrix}
      \frac{g_{k-1} \tr \Psi^{2r}(A)Ag_{k-1}}{\|\Psi^r(A)g_{k-1}\|^2} & -\frac{g_k \tr \Psi(A) g_k}{\alpha_{k-1}\|\Psi^r(A)g_{k-1}\|\|\Psi^{1-r}(A)g_k\|} \\
      -\frac{g_k \tr \Psi(A) g_k}{\alpha_{k-1}\|\Psi^r(A)g_{k-1}\|\|\Psi^{1-r}(A)g_k\|}  & \frac{g_k \tr \Psi^{2(1-r)}(A)Ag_k}{\|\Psi^{1-r}(A)g_k\|^2} \\
    \end{pmatrix}
\end{equation}
and $\alpha_{k-1}$ is given by \eqref{glstep1}.
Clearly, $\alpha_k^Y$ in \eqref{syv} can be obtained by
by setting $\Psi(A)=I$ in \eqref{eqhesk}.
In addition, by \eqref{newn} we have that
\begin{equation}\label{newbd}
\frac{1}{H_{11}^k+H_{22}^k}
  \leq\tilde{\alpha}_k
  \leq\frac{1}{\max\{H_{11}^k,H_{22}^k\}}.
\end{equation}
%Thus, for all $r$ such that the minimum order of the series of $\Psi^{1-r}(A)$ is positive, we have that $%\tilde{\alpha}_k\leq\alpha_k^{SD}$,
%, the new stepsize $\tilde{\alpha}_k$ will decrease the objective value.

%To analyze the property of $\tilde{\alpha}_k$, we first consider the asymptotic behavior of $\alpha_{k}^{SD}$. The next corollary can be obtained from Theorem \ref{th1}.
%
%
%\begin{corollary}\label{cor2}
%Under the conditions of Theorem \ref{th1}, the following limits hold,
%\begin{equation}\label{sd1}
%  \lim_{k\rightarrow\infty}\alpha_{2k}^{SD}=\frac{1+c^2}{\lambda_1(1+c^2\kappa)}
%\end{equation}
%and
%\begin{equation}\label{sd2}
%  \lim_{k\rightarrow\infty}\alpha_{2k+1}^{SD}=
%  \frac{\Psi^2(\lambda_{1})+c^2\Psi^2(\lambda_{n})}{\lambda_1(\kappa\Psi^2(\lambda_{1})+c^2\Psi^2(\lambda_{n}))}.
%\end{equation}
%%where $c$ is the same constant as in Theorem \ref{th1}.
%\end{corollary}

The next theorem shows that the stepsize $\tilde{\alpha}_k$ enjoys desirable spectral property.

\begin{theorem}
\label{thtilalp}
Suppose that the conditions of Theorem \ref{th1} hold. Let $\{x_k\}$ be the iterations
generated by any gradient method in \eqref{glstep1} to solve problem \eqref{qudpro}. Then
\begin{equation}\label{tidlbdn}
  \lim_{k\rightarrow\infty}\tilde{\alpha}_k=\frac{1}{\lambda_n}.
\end{equation}
\end{theorem}
\begin{proof}
  It follows from \eqref{mu2k} and \eqref{mu2k1} of Theorem \ref{th1} that
  \begin{align*}
  \lim_{k\rightarrow\infty}H_{11}^k
  %\frac{g_{k-1} \tr \Psi^{2r}(A)Ag_{k-1}}{\|\Psi^r(A)g_{k-1}\|^2}
  &=\lim_{k\rightarrow\infty}\frac{g_{k-1} \tr \Psi^{2r}(A)Ag_{k-1}}{\|g_{k-1}\|^2}
  \frac{\|g_{k-1}\|^2}{\|\Psi^r(A)g_{k-1}\|^2}\\
  &=\frac{\lambda_1(c^2\Psi^{2r}(\lambda_{1})\Psi^2(\lambda_{n})+\kappa\Psi^{2r}(\lambda_{n})\Psi^2(\lambda_{1}))}
  {c^2\Psi^{2r}(\lambda_{1})\Psi^2(\lambda_{n})+\Psi^{2r}(\lambda_{n})\Psi^2(\lambda_{1})}
\end{align*}
and
\begin{align*}
  \lim_{k\rightarrow\infty}H_{22}^k&=
  \frac{g_k \tr \Psi^{2(1-r)}(A)Ag_k}{\|g_k\|^2}
  \frac{\|g_k\|^2}{\|\Psi^{1-r}(A)g_k\|^2}  \\
  &=\frac{\lambda_1(\Psi^{2(1-r)}(\lambda_{1})+\kappa c^2\Psi^{2(1-r)}(\lambda_{n}))}{\Psi^{2(1-r)}(\lambda_{1})+ c^2\Psi^{2(1-r)}(\lambda_{n})}\\
    &=\frac{\lambda_1(\Psi^2(\lambda_{1})\Psi^{2r}(\lambda_{n})+\kappa c^2\Psi^2(\lambda_{n})\Psi^{2r}(\lambda_{1}))}
    {\Psi^2(\lambda_{1})\Psi^{2r}(\lambda_{n})+ c^2\Psi^2(\lambda_{n})\Psi^{2r}(\lambda_{1})},
\end{align*}
which give
  \begin{equation}\label{mgsd}
  \lim_{k\rightarrow\infty}(H_{11}^k+H_{22}^k)
  =\lambda_1(\kappa+1)
\end{equation}
and
  \begin{equation}\label{mgsd2}
  \lim_{k\rightarrow\infty}(H_{11}^k-H_{22}^k)
  =\frac{\lambda_1(\kappa-1)(\Psi^2(\lambda_{1})\Psi^{2r}(\lambda_{n})- c^2\Psi^2(\lambda_{n})\Psi^{2r}(\lambda_{1}))}
    {\Psi^2(\lambda_{1})\Psi^{2r}(\lambda_{n})+ c^2\Psi^2(\lambda_{n})\Psi^{2r}(\lambda_{1})}.
\end{equation}
Then, by the definition of $\alpha_k$, we have
\[
  g_{k} \tr \Psi(A)g_{k}=-\alpha_{k-1}g_{k-1} \tr \Psi(A)Ag_{k-1}
  +\alpha_{k-1}^2g_{k-1} \tr \Psi(A)A^2g_{k-1},
\]
which together with \eqref{mu2k1} in Theorem \ref{th1} and \eqref{salp2} in Corollary \ref{cor1} yields that
\begin{align*}
&\lim_{k\rightarrow\infty}(H_{12}^k)^2\\
   =&\lim_{k\rightarrow\infty}\frac{g_{k} \tr \Psi(A)g_{k}}{\alpha_{k-1}^2\|\Psi^r(A)g_{k-1}\|^2}
   \frac{g_k \tr \Psi(A) g_k}{\|\Psi^{1-r}(A)g_k\|^2}\\
= &\lim_{k\rightarrow\infty}
   \left(-\frac{1}{\alpha_{k-1}}\frac{g_{k-1} \tr \Psi(A)Ag_{k-1}}{\|\Psi^r(A)g_{k-1}\|^2}
  +\frac{g_{k-1} \tr \Psi(A)A^2g_{k-1}}{\|\Psi^r(A)g_{k-1}\|^2}\right)
   \frac{g_k \tr \Psi(A) g_k}{\|\Psi^{1-r}(A)g_k\|^2}\\
 % =&\left(-\frac{\lambda_1^2(\Psi(\lambda_{1})+c^2\kappa\Psi(\lambda_{n}))^2}
%  {(1+c^2)(\Psi(\lambda_{1})+c^2\Psi(\lambda_{n}))}
%  +
%  \frac{\lambda_1^2(\Psi(\lambda_{1})+c^2\kappa^2\Psi(\lambda_{n}))}
%  {(1+c^2)}\right)\\
  =&
\Bigg[-\frac{\lambda_{1}(\kappa\Psi(\lambda_{1})+c^2\Psi(\lambda_{n}))}{\Psi(\lambda_{1})+c^2\Psi(\lambda_{n})}
  \frac{\lambda_1(c^2\Psi(\lambda_{1})\Psi^2(\lambda_{n})+\kappa\Psi(\lambda_{n})\Psi^2(\lambda_{1}))}
  {c^2\Psi^{2r}(\lambda_{1})\Psi^2(\lambda_{n})+\Psi^{2r}(\lambda_{n})\Psi^2(\lambda_{1})}
  +\\
  &
\frac{\lambda_1^2(c^2\Psi(\lambda_{1})\Psi^2(\lambda_{n})+\kappa^2\Psi(\lambda_{n})\Psi^2(\lambda_{1}))}
  {c^2\Psi^{2r}(\lambda_{1})\Psi^2(\lambda_{n})+\Psi^{2r}(\lambda_{n})\Psi^2(\lambda_{1})}
  \Bigg]
  \frac{(\Psi(\lambda_{1})+ c^2\Psi(\lambda_{n}))\Psi^{2r}(\lambda_{1})\Psi^{2r}(\lambda_{n})}{\Psi^2(\lambda_{1})\Psi^{2r}(\lambda_{n})+ c^2\Psi^2(\lambda_{n})\Psi^{2r}(\lambda_{1})}
%  \frac{\Psi(\lambda_{1})+ c^2\Psi(\lambda_{n})}{\Psi^{2(1-r)}(\lambda_{1})+ c^2\Psi^{2(1-r)}(\lambda_{n})}
  \\
 % &
%  =\Big[-\frac{1}{(c^2\Psi(\lambda_{1})\Psi^2(\lambda_{n})+\Psi(\lambda_{n})\Psi^2(\lambda_{1}))}
%  \frac{\lambda_1^2(c^2\Psi(\lambda_{1})\Psi^2(\lambda_{n})+\kappa\Psi(\lambda_{n})\Psi^2(\lambda_{1}))^2}
%  {c^2\Psi^{2r}(\lambda_{1})\Psi^2(\lambda_{n})+\Psi^{2r}(\lambda_{n})\Psi^2(\lambda_{1})}\\
%  &
%  +\frac{\lambda_1^2(c^2\Psi(\lambda_{1})\Psi^2(\lambda_{n})+\kappa^2\Psi(\lambda_{n})\Psi^2(\lambda_{1}))}
%  {c^2\Psi^{2r}(\lambda_{1})\Psi^2(\lambda_{n})+\Psi^{2r}(\lambda_{n})\Psi^2(\lambda_{1})}
%  \Big]
%  \frac{(\Psi(\lambda_{1})+ c^2\Psi(\lambda_{n}))\Psi^{2r}(\lambda_{1})\Psi^{2r}(\lambda_{n})}{\Psi^2(\lambda_{1})\Psi^{2r}(\lambda_{n})+ c^2\Psi^2(\lambda_{n})\Psi^{2r}(\lambda_{1})}
%  \\
  =& \frac{\lambda_1^2c^2(\kappa-1)^2\Psi^{2+2v}(\lambda_{1})\Psi^{2+2v}(\lambda_{n})}
  {(\Psi^2(\lambda_{1})\Psi^{2r}(\lambda_{n})+ c^2\Psi^2(\lambda_{n})\Psi^{2r}(\lambda_{1}))^2}.
\end{align*}
Then, from the above equality and \eqref{mgsd2}, we obtain that
  \begin{equation}\label{mgsd3}
  \lim_{k\rightarrow\infty}\sqrt{(H_{11}^k-H_{22}^k)^2+4(H_{12}^k)^2}
  =\lambda_1(\kappa-1).
\end{equation}
Combining \eqref{mgsd} and \eqref{mgsd3}, we have that
\[
  \lim_{k\rightarrow\infty}\tilde{\alpha}_k=
  \frac{2}{\lambda_1(\kappa+1)+\lambda_1(\kappa-1)}=\frac{1}{\lambda_n}.
\]
This completes the proof.
\end{proof}

\begin{remark}
When $r=1$, we have from \eqref{newbd} that $\tilde{\alpha}_k \le 1/H_{22}^k = \alpha_k^{SD}$.
Hence,  using this stepsize $\tilde{\alpha}_k$ will give a monotone  gradient method.
Theorem \ref{thtilalp} indicates that the general $\tilde{\alpha}_k$ will have
the asymptotic spectral property \eqref{tidlbdn}, and hence  will be asymptotically
be smaller than $\alpha_k^{SD}$ independent of $r$.
But a proper choice $r$ will facilitate the calculation of
 $\tilde{\alpha}_k$. This will be more clear in the next section.
\end{remark}

Using the similar arguments, we can also show the larger stepsize derived in
 subsection \ref{subsecns} converges to $1/\lambda_1$.
\begin{theorem}
Let
\[
    \bar{\alpha}_k=
  \frac{2}{(H_{11}^k+H_{22}^k)-\sqrt{(H_{11}^k-H_{22}^k)^2+4(H_{12}^k)^2}}.
\]
Under the conditions of Theorem \ref{thtilalp}, we have
\[
  \lim_{k\rightarrow\infty}\bar{\alpha}_k=\frac{1}{\lambda_1}.
\]
\end{theorem}

\begin{figure}[h]
  \centering
  \includegraphics[width=0.65\textwidth,height=0.42\textwidth]{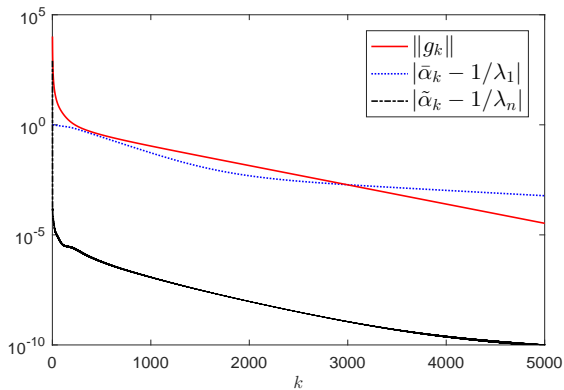}\\
  \caption{Problem \eqref{tp1} with $n=1,000$: convergence history of the sequences $\{\tilde{\alpha}_k\}$ and $\{\bar{\alpha}_k\}$ for the first 5,000 iterations of the gradient method \eqref{glstep1} with $\Psi(A)=A$ (i.e., the MG method).}\label{appstep1}
\end{figure}

To present an intuitive illustration of the asymptotic behaviors of $\tilde{\alpha}_k$ and $\bar{\alpha}_k$, we applied the gradient method
 \eqref{glstep1} with $\Psi(A)=A$ (i.e., the MG method) to minimize
the quadratic function \eqref{qudpro} with
\begin{equation}\label{tp1}
  A=\textrm{diag}\{a_1,a_2,\ldots,a_n\} \quad \mbox{and} \quad b=0,
\end{equation}
where $a_1=1$, $a_n=n$ and $a_i$ is randomly generated between 1 and $n$ for $i=2,\ldots,n-1$. From Figure \ref{appstep1}, we can see that $\tilde{\alpha}_k$ approximates $1/\lambda_n$ with satisfactory accuracy in a few iterations. However, $\bar{\alpha}_{k}$ converges to $1/\lambda_1$ even slower than the decreasing of gradient norm. This, to some extent, explains the reason why we prefer $\tilde{\alpha}_k$ to the
 short stepsize.

\subsection{A periodic gradient method}\label{secmethod}

A method alternately using $\alpha_k$ in \eqref{glstep1} and $\tilde{\alpha}_k$ to minimize a
$2$-dimensional quadratic function  will monotonically decrease the objective value, and
terminates in $3$ iterations.
 However, for minimizing a general $n$-dimensional quadratic function, this alternating scheme may not be efficient for the purpose of vanishing the component $\mu_k^{(n)}$.
 One possible reason is that, as shown in Figure \ref{appstep1}, it needs tens of iterations before
  $\tilde{\alpha}_k$ being a good approximation of $1/\lambda_n$ with satisfactory accuracy.
In what follows, by incorporating the BB method, we develop an efficient periodic gradient method
using $\tilde{\alpha}_k$.

\begin{figure}[tbhp]
\centering
\subfloat{\label{fig:a}\includegraphics[width=0.47\textwidth,height=0.37\textwidth]{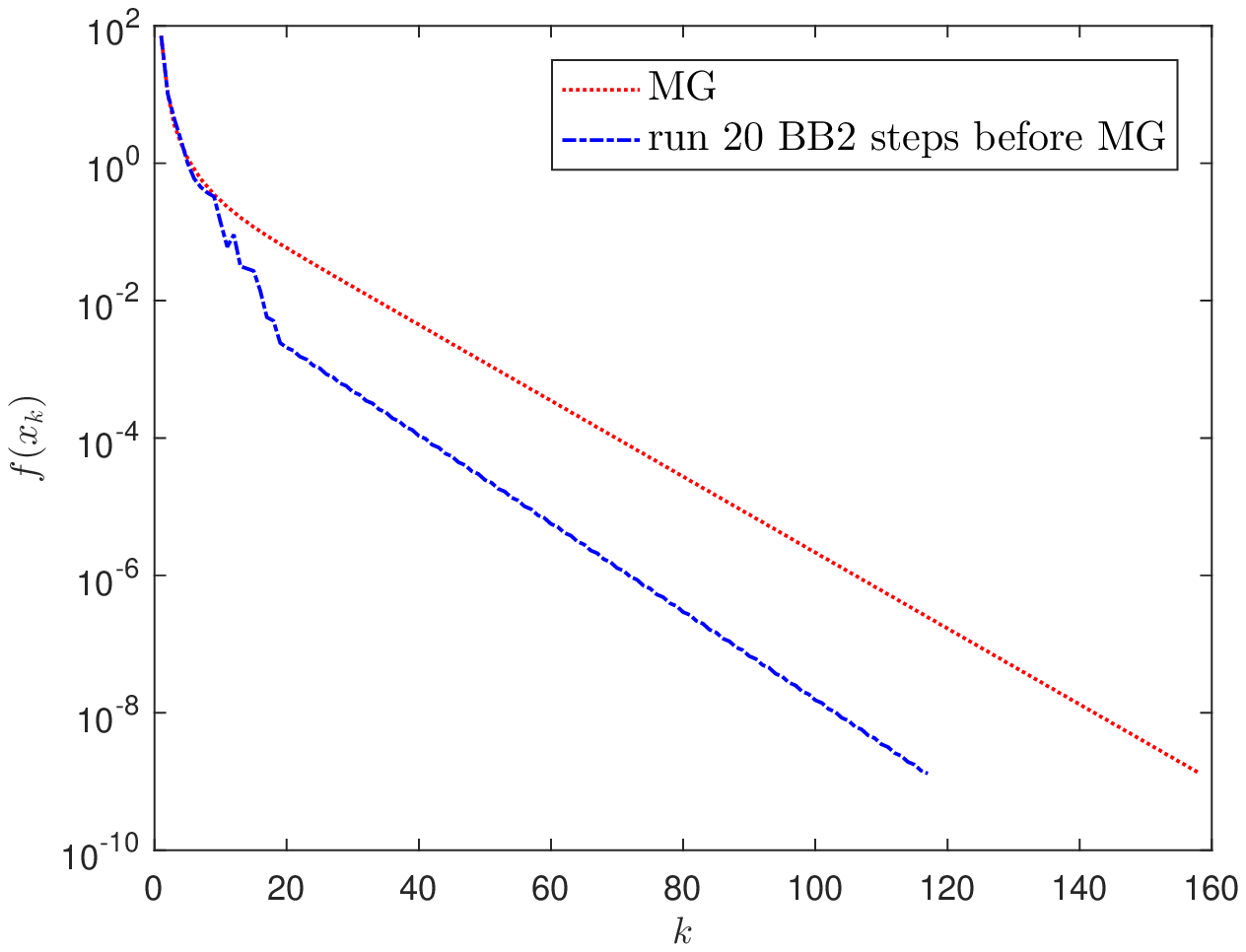}}
\subfloat{\label{fig:b}\includegraphics[width=0.47\textwidth,height=0.37\textwidth]{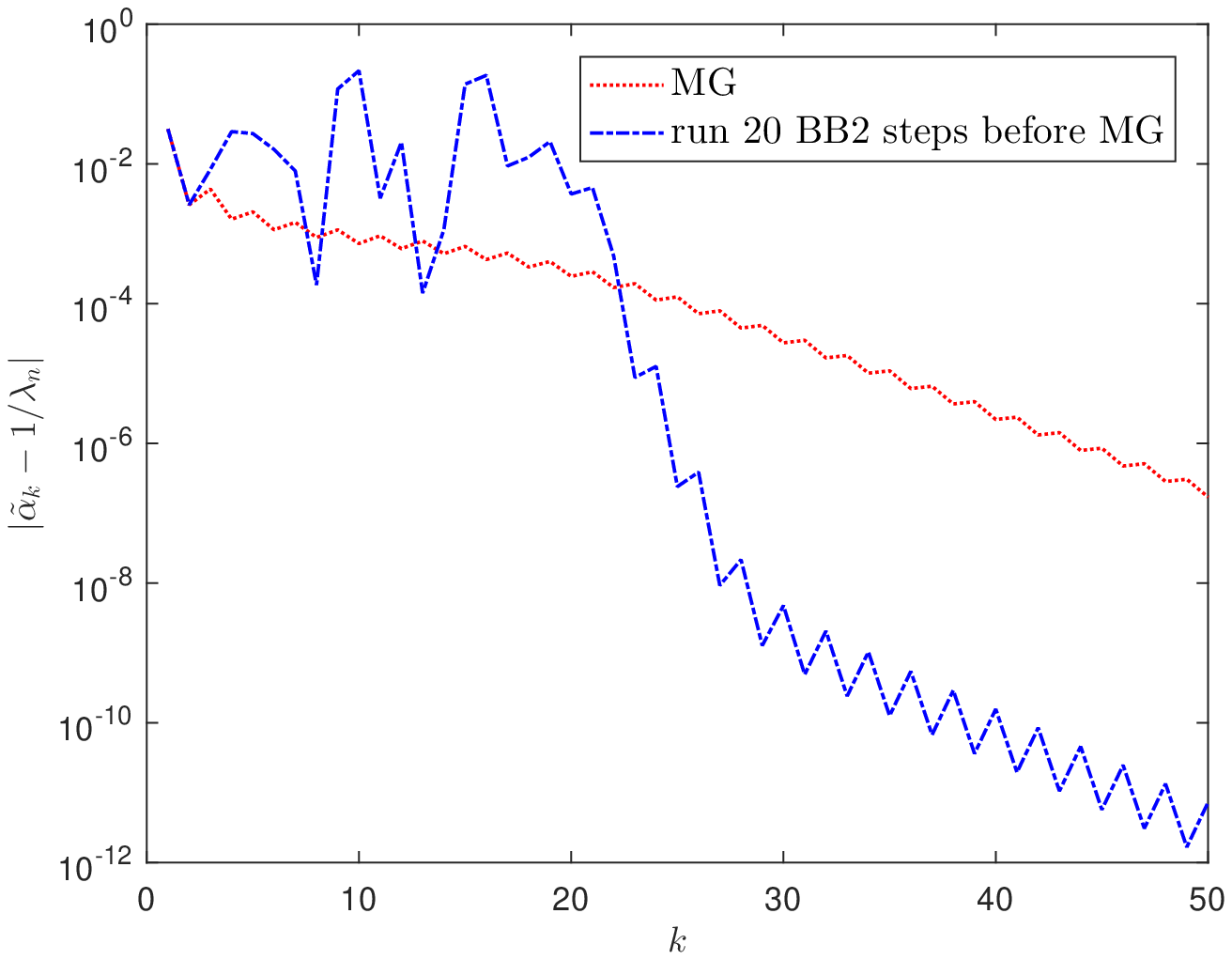}}
  \caption{Problem \eqref{exfkrate} with $n=10$: convergence history of objective values and stepsizes.}\label{BBMG}
\end{figure}

\begin{figure}[h]
  \centering
  \includegraphics[width=0.6\textwidth,height=0.42\textwidth]{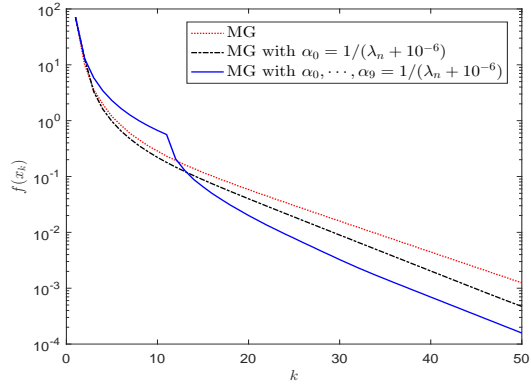}\\
  \caption{Problem \eqref{exfkrate} with $n=10$: the MG method (i.e., $\Psi(A)=A$) with different stepsizes.}\label{BBMGstep}
\end{figure}

Figure \ref{BBMG} illustrates a comparison of the gradient method \eqref{glstep1} using $\Psi(A)=A$
(i.e., the MG method) with a method using 20 BB2 steps first and then MG steps on solving
 problem \eqref{exfkrate}.  We can see that by  using some BB2 steps,
 the modified MG method is accelerated and the stepsize $\tilde{\alpha}_k$ will approximate $1/\lambda_n$
  with a better accuracy. Thus, our  method will run some BB steps first.
Now, we investigate the affect of reusing a short stepsize on the performance of the gradient method \eqref{glstep1}. Suppose that we have a good approximation of $1/\lambda_n$, say $\alpha=\frac{1}{\lambda_n+10^{-6}}$. We compare MG method with its two variants by applying
(i) $\alpha_0=\alpha$ or (ii) $\alpha_0=\ldots=\alpha_9=\alpha$ before using the MG stepsize.
Figure \ref{BBMGstep} shows that reusing $\alpha$ will accelerate the MG method.
Hence, we prefer to reuse $\tilde{\alpha}_k$ for some consecutive steps
when $\tilde{\alpha}_k$ is a good approximation of $1/\lambda_n$. Finally,
our new method is summarized in Algorithm \ref{al1}, which periodically applies the BB stepsize,
 $\alpha_k$ in \eqref{glstep1} and $\tilde{\alpha}_k$.
 The $R$-linear global convergence of Algorithm \ref{al1} for solving \eqref{qudpro} can
 be established by showing that it satisfies the property in \cite{dai2003alternate}, see Theorem 3 of \cite{dhl2018} for example.

%\begin{figure}[h]
%  \centering
%  \subfigure{
% \includegraphics[width=0.47\textwidth,height=0.38\textwidth]{BB2_MG_largeg_obj.pdf}}
%\subfigure{
%    \includegraphics[width=0.47\textwidth,height=0.38\textwidth]{BB2_MG_largeg.pdf}}
%  \caption{Problem \eqref{exfkrate} with $n=10$: convergence history of objective values (left) and stepsizes (right).}\label{BBMG}
%\end{figure}

\begin{algorithm}[ht!b]
\caption{Periodic gradient method}\label{al1}
\begin{algorithmic}
\STATE{Choose an initial point $x_{0}\in \mathbb{R}^{n}$, initial stepsize $\alpha_{0}$, positive integers $K_b,K_m,K_s$, and termination tolerance $\epsilon>0$.}

\STATE{Take one gradient step with $\alpha_{0}$}

\WHILE{$\|g_k\|>\epsilon$}
\STATE{Take $K_b$ BB steps}
\STATE{Take $K_m$ gradient steps with $\alpha_k$ in \eqref{glstep1} }
\STATE{Take $K_s$ short steps with $\tilde{\alpha}_t$, where $\tilde{\alpha}_t$ is the first stepsize after $\alpha_k$-steps}
\ENDWHILE
%\RETURN $T$
\end{algorithmic}
\end{algorithm}

%%\SetAlgoSkip{bigskip}
%\DontPrintSemicolon
%\begin{algorithm}[ht!b]
%\caption{Periodic gradient method}\label{al1}
%Choose an initial point $x_{0}\in \mathbb{R}^{n}$, initial stepsize $\alpha_{0}$, positive integers $K_b,K_m,K_s$, and termination tolerance $\epsilon>0$.\\
%
%Take one gradient step with $\alpha_{0}$.
%
%\Repeat{$\|g_k\|\leq\epsilon$.}{
%
%Take $K_b$ BB steps.
%
%Take $K_m$ gradient steps with $\alpha_k$ in \eqref{glstep1}.
%
%Take $K_s$ short steps with $\tilde{\alpha}_t$, where $\tilde{\alpha}_t$ is the first stepsize after $\alpha_k$-steps.
%}
%\end{algorithm}

\begin{remark}
The BB steps in Algorithm \ref{al1} can either employ the BB1 or BB2 stepsize in \eqref{sBB}.
  The idea of using short stepsizes to eliminate the component $\mu_k^{(n)}$ has been investigated in \cite{de2014efficient,de2013spectral,gonzaga2016steepest}. However, these methods are based on the SD method, that is, occasionally applying short steps during the iterates of the SD method. One exception is given by \cite{huang2019gradient}, where a method is developed by employing new stepsizes during the iterates of the AOPT method.
But our method periodically uses three different stepsizes: 
 the nonmonotone BB method, the gradient method \eqref{glstep1} and the new stepsize $\tilde{\alpha}_k$.
\end{remark}

\section{Numerical experiments}\label{secnum}
In this section, we present numerical comparisons of Algorithm \ref{al1} and the following methods: BB with $\alpha_k^{BB1}$ \cite{barzilai1988two}, Dai-Yuan (DY) \cite{dai2005analysis}, ABBmin2 \cite{frassoldati2008new}, and SDC \cite{de2014efficient}.
%Since the SDC method performs better than its monotone counterpart for most problems, we only run SDC.

Notice that the stepsize rule for Algorithm \ref{al1} can be written as
\begin{equation}\label{sal1}
\alpha_{k}=
\begin{cases}
\alpha_{k}^{BB},& \text{if $\textrm{mod}(k,K_b+K_m+K_s)<K_b$}; \\
\alpha_{k}(\Psi(A)),& \text{if $K_b\leq\textrm{mod}(k,K_b+K_m+K_s)<K_b+K_m$}; \\
\tilde{\alpha}_{k}(\Psi(A)),& \text{if $\textrm{mod}(k,K_b+K_m+K_s)=K_b+K_m$}; \\
\alpha_{k-1},& \text{otherwise},
\end{cases}
\end{equation}
where $\alpha_k^{BB}$ can either be $\alpha_k^{BB1}$ or $\alpha_k^{BB2}$, $\alpha_{k}(\Psi(A))$ and $\tilde{\alpha}_{k}(\Psi(A))$ are the stepsizes given by \eqref{glstep1} and \eqref{newn}, respectively. We tested the following four variants of Algorithm \ref{al1} using combinations of the two BB stepsizes and $\Psi(A)=I$ or $A$:
\begin{itemize}
  \item BB1SD: $\alpha_k^{BB1}$ and $\Psi(A)=I$ in \eqref{sal1}
  \item BB1MG: $\alpha_k^{BB1}$ and $\Psi(A)=A$ in \eqref{sal1}
  \item BB2SD: $\alpha_k^{BB2}$ and $\Psi(A)=I$ in \eqref{sal1}
  \item BB2MG: $\alpha_k^{BB2}$ and $\Psi(A)=A$ in \eqref{sal1}
\end{itemize}

Now we derive a formula for the case $\Psi(A)=A$, i.e., $\alpha_{k}(\Psi(A))=\alpha_{k}^{MG}$. If we set $r=0$, by \eqref{newn}, we have
\begin{equation}\label{newmg1}
  \tilde{\alpha}_k=\frac{2}{\left(\frac{1}{\alpha_{k-1}^{SD}}+ \frac{g_k \tr A^3g_k}{g_k \tr A^2g_k}\right)+
  \sqrt{\left(\frac{1}{\alpha_{k-1}^{SD}}- \frac{g_k \tr A^3g_k}{g_k \tr A^2g_k}\right)^2+
  \frac{4(g_k \tr Ag_k)^2}{(\alpha_{k-1}^{MG})^2\|g_{k-1}\|^2g_k \tr A^2g_k}}},
\end{equation}
which is expensive to compute directly.
%By using
%\begin{equation*}
%  \frac{g_{k} \tr Ag_{k}}{g_{k-1} \tr Ag_{k-1}}=
%  \frac{\alpha_{k-1}^{SD}-\alpha_{k-1}^{MG}}{\alpha_{k}^{SD}}
%\end{equation*}
%and
%\begin{equation*}
%  \frac{g_{k-1} \tr A^3g_{k-1}}{g_{k-1} \tr A^2g_{k-1}}=
%  \frac{1+\frac{\alpha_{k-1}^{SD}-\alpha_{k-1}^{MG}}{\alpha_{k}^{SD}}}{\alpha_{k-1}^{MG}},
%\end{equation*}
%we obtain the following formula for $\tilde{\alpha}_k$:
%\begin{equation}\label{newmg2}
%  \tilde{\alpha}_k=\frac{2\alpha_{k-1}^{MG}\alpha_{k}^{SD}}{(\alpha_{k-1}^{SD}+ \alpha_{k}^{SD})+
%  \sqrt{(\alpha_{k-1}^{SD}+\alpha_{k}^{SD})^2-4\alpha_{k-1}^{MG}\alpha_{k}^{SD}}}.
%\end{equation}
%This formula can be computed as cheap as $\alpha_{k}^{SD}$, i.e., by only one matrix-vector product in each iteration. In fact, from
%\begin{equation*}
%  g_{k}=g_{k-1}-\alpha_{k-1}^{MG}Ag_{k-1},
%\end{equation*}
%we have that
%\begin{equation*}
%  \alpha_{k-1}^{SD}=\frac{g_{k-1} \tr g_{k-1}}{g_{k-1} \tr Ag_{k-1}}
%  =\alpha_{k-1}^{MG}\frac{g_{k-1} \tr g_{k-1}}{g_{k-1} \tr (g_{k-1}-g_{k})}.
%\end{equation*}
%Notice that $g_{k} \tr Ag_{k-1}=0$. We further have
%\begin{equation*}
%  g_{k} \tr g_{k}=g_{k-1} \tr g_{k}-\alpha_{k-1}^{MG}g_{k} \tr Ag_{k-1}
%  =g_{k-1} \tr g_{k}.
%\end{equation*}
%Thus, we only need to compute $g_{k} \tr Ag_{k}$ for $\alpha_{k}^{SD}$.
However, if we set $r=1/2$, we get
\begin{equation}\label{newsmg}
  \tilde{\alpha}_{k}=
  \frac{2}{\frac{1}{\alpha_{k-1}^{MG}}+\frac{1}{\alpha_k^{MG}}+
  \sqrt{\left(\frac{1}{\alpha_{k-1}^{MG}}-\frac{1}{\alpha_k^{MG}}\right)^2
  +\frac{4g_k \tr A g_k}{(\alpha_{k-1}^{MG})^2g_{k-1} \tr A g_{k-1}}}}.
\end{equation}
This formula can be computed without additional cost because $g_{k-1} \tr A g_{k-1}$ and $g_k \tr A g_k$ have been obtained when computing the stepsizes $\alpha_{k-1}^{MG}$ and  $\alpha_k^{MG}$.

All the methods in consideration were implemented in Matlab (v.9.0-R2016a) and carried out on a PC with an Intel Core i7, 2.9 GHz processor and 8 GB of RAM running Windows 10 system. We stopped the algorithm if the number of iteration exceeds 20,000 or the gradient norm reduces by a factor of $\epsilon$.

We randomly generated quadratic problems \eqref{eqpro} proposed in \cite{dhl2018}, where $A=QVQ \tr $  with
\[
  Q=(I-2w_3w_3 \tr )(I-2w_2w_2 \tr )(I-2w_1w_1 \tr ),
\]
 $w_1$, $w_2$, and $w_3$ are unitary random vectors, and $V=diag(v_1,\ldots,v_n)$ is a diagonal matrix where $v_1=1$, $v_n=\kappa$, and $v_j$, $j=2,\ldots,n-1$, are randomly generated between 1 and $\kappa$ by the \emph{rand} function in Matlab.
We tested seven sets of different distributions of $v_j$  as shown in Table \ref{tbspe} with different values of the condition number $\kappa$ and tolerance $\epsilon$. In particular, $\kappa$ were set to $10^4, 10^5, 10^6$ and $\epsilon$ were set to $10^{-6}, 10^{-9}, 10^{-12}$. For each value of $\kappa$ or $\epsilon$, 10 instances were generated and there are totally 630 instances. For each instance, the entries of $b$ were randomly generated in $[-10,10]$ and
$e=(1,\ldots,1) \tr$ was used as the starting point.

\begin{table}[ht!b]
{\footnotesize
\caption{Distributions of $v_j$.}\label{tbspe}
\begin{center}
\begin{tabular}{|c|c|}
 \hline
 \multirow{1}{*}{Set} &\multicolumn{1}{c|}{Spectrum} \\
\hline
\multirow{1}{*}{1} &$\{v_2,\ldots,v_{n-1}\}\subset(1,\kappa)$	\\
\hline
 \multirow{2}{*}{2}
&$\{v_2,\ldots,v_{n/5}\}\subset(1,100)$	\\
&$\{v_{n/5+1},\ldots,v_{n-1}\}\subset(\frac{\kappa}{2},\kappa)$	\\
 \hline

\multirow{2}{*}{3}
&$\{v_2,\ldots,v_{n/2}\}\subset(1,100)$	\\
&$\{v_{n/2+1},\ldots,v_{n-1}\}\subset(\frac{\kappa}{2},\kappa)$	\\
 \hline
\multirow{2}{*}{4}
&$\{v_2,\ldots,v_{4n/5}\}\subset(1,100)$	\\
&$\{v_{4n/5+1},\ldots,v_{n-1}\}\subset(\frac{\kappa}{2},\kappa)$	\\
 \hline
\multirow{3}{*}{5}
&$\{v_2,\ldots,v_{n/5}\}\subset(1,100)$	\\
&$\{v_{n/5+1},\ldots,v_{4n/5}\}\subset(100,\frac{\kappa}{2})$	\\
&$\{v_{4n/5+1},\ldots,v_{n-1}\}\subset(\frac{\kappa}{2},\kappa)$	\\
 \hline
 \multirow{2}{*}{6}
&$\{v_2,\ldots,v_{10}\}\subset(1,100)$	\\
&$\{v_{11},\ldots,v_{n-1}\}\subset(\frac{\kappa}{2},\kappa)$	\\
 \hline
 \multirow{2}{*}{7}
&$\{v_2,\ldots,v_{n-10}\}\subset(1,100)$	\\
&$\{v_{n-9},\ldots,v_{n-1}\}\subset(\frac{\kappa}{2},\kappa)$	\\
 \hline
\end{tabular}
\end{center}
}
\end{table}

%\begin{table}[ht!b]
%%\setlength{\tabcolsep}{0.3ex}
%\begin{center}
%%\begin{scriptsize}
%\caption{Distributions of $v_j$.}\label{tbspe}
%%\noindent\makebox[\textwidth]{
%\begin{tabular}{|c|c|c|c|}
% \hline
% \multirow{1}{*}{Set} &\multicolumn{1}{c|}{Spectrum}
% &\multirow{1}{*}{Set} &\multicolumn{1}{c|}{Spectrum}\\
%\hline
%\multirow{3}{*}{1} & 	&\multirow{3}{*}{5}
%&$\{v_2,\ldots,v_{n/5}\}\subset(1,100)$	\\
%&$\{v_2,\ldots,v_{n-1}\}\subset(1,\kappa)$ &
%&$\{v_{n/5+1},\ldots,v_{4n/5}\}\subset(100,\frac{\kappa}{2})$	\\
%& 	&
%&$\{v_{4n/5+1},\ldots,v_{n-1}\}\subset(\frac{\kappa}{2},\kappa)$	\\
%\hline
%
% \multirow{2}{*}{2}
%&$\{v_2,\ldots,v_{n/5}\}\subset(1,100)$	& \multirow{2}{*}{6}
%&$\{v_2,\ldots,v_{10}\}\subset(1,100)$	\\
%&$\{v_{n/5+1},\ldots,v_{n-1}\}\subset(\frac{\kappa}{2},\kappa)$	
%&
%&$\{v_{11},\ldots,v_{n-1}\}\subset(\frac{\kappa}{2},\kappa)$	\\
% \hline
%
%\multirow{2}{*}{3}
%&$\{v_2,\ldots,v_{n/2}\}\subset(1,100)$	 &\multirow{2}{*}{7}
%&$\{v_2,\ldots,v_{n-10}\}\subset(1,100)$	\\
%
%&$\{v_{n/2+1},\ldots,v_{n-1}\}\subset(\frac{\kappa}{2},\kappa)$	
% &
%&$\{v_{n-9},\ldots,v_{n-1}\}\subset(\frac{\kappa}{2},\kappa)$	\\
% \hline
%\multirow{2}{*}{4}
%&$\{v_2,\ldots,v_{4n/5}\}\subset(1,100)$	& &\\
%&$\{v_{4n/5+1},\ldots,v_{n-1}\}\subset(\frac{\kappa}{2},\kappa)$
%& &\\
% \hline
%\end{tabular}
%%\end{scriptsize}
%\end{center}
%\end{table}

The parameter $K_b$ for Algorithm \ref{al1} was set to 100 for the first and fifth sets and 30 for other sets. Other two parameters $K_m$ and $K_s$ were selected from $\{9,13,15\}$. As in \cite{frassoldati2008new}, the parameter $\tau$ of the ABBmin2 method was set to 0.9 for all instances. The parameter pair $(h,s)$ used for the SDC method was set to $(8,6)$, which is more efficient than other choices for this test.

Table \ref{tbrandpBB1SD} shows the averaged number of iterations of BB1SD and other four compared methods for the seven sets of problems listed in Table \ref{tbspe}. We can see that, for the first problem set, our BB1SD method performs much better than the BB, DY and SDC methods, although the ABBmin2 method seems surprisingly efficient among the compared methods. For the second to the last problem sets, our method with different settings performs better than the BB, DY, ABBmin2 and SDC methods. Moreover, for all the settings and different tolerance levels, our method outperforms all the compared four methods in terms of total number of iterations.

Tables \ref{tbrandpBB1MG}, \ref{tbrandpBB2SD} and \ref{tbrandpBB2MG} present the averaged number of iterations of BB1MG, BB2SD and BB2MG, respectively. For comparison purposes, the results of the BB, DY, ABBmin 2 and SDC methods are also listed in those tables. As compared with the BB, DY, ABBmin 2 and SDC methods, similar results to those in Table \ref{tbrandpBB1SD} can be seen from these three tables. For the comparison of BB1SD and BB1MG, we can see from Tables \ref{tbrandpBB1SD} and \ref{tbrandpBB1MG} that BB1MG is slightly better than BB1SD for the second to fourth, sixth, and the last problem sets. In addition, BB1MG is comparable to BB1SD for the first and the fifth problem sets. The results in Tables \ref{tbrandpBB2SD} and \ref{tbrandpBB2MG} do not show much difference between BB2SD and BB2MG. In general, BB1MG performs slightly better than BB1SD, BB2SD and BB2MG for most of the problem sets.

\begin{figure}[tbhp]
\centering
\subfloat{\label{fig:2a}\includegraphics[width=0.47\textwidth,height=0.37\textwidth]{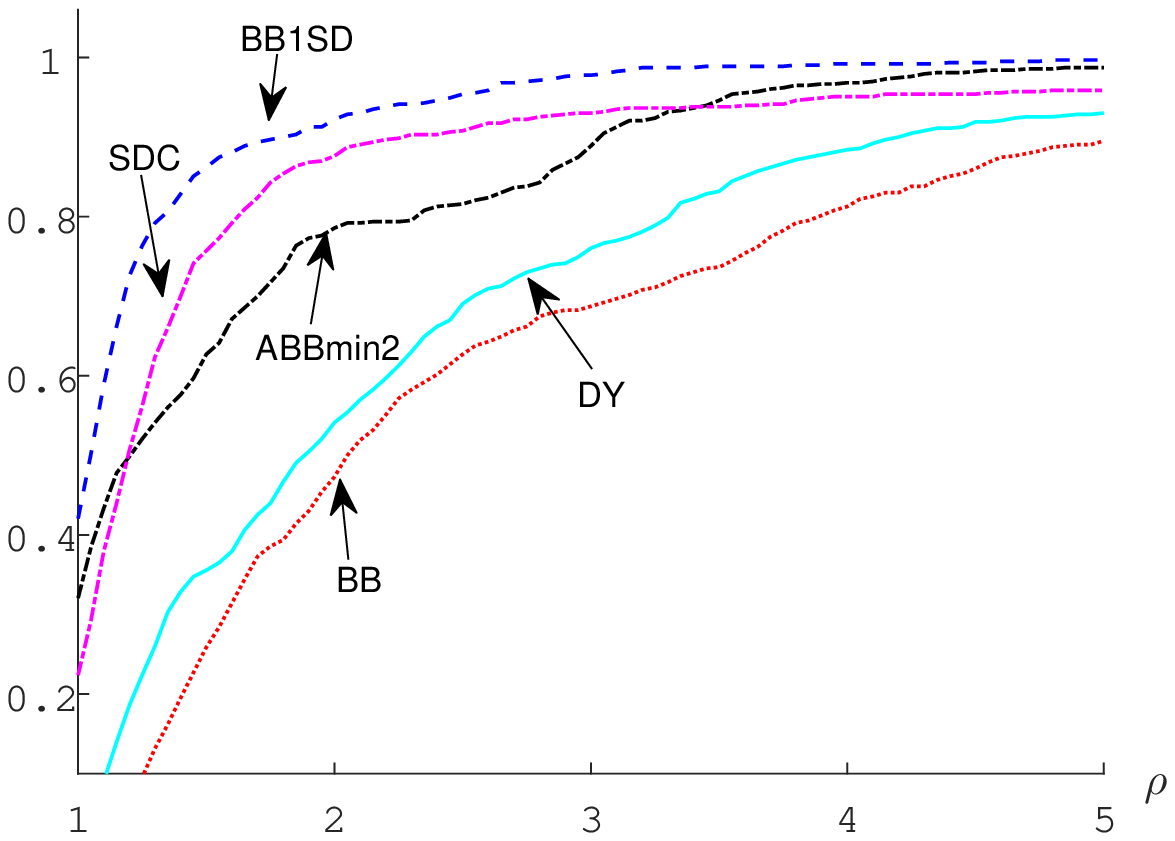}}
\subfloat{\label{fig:2b}\includegraphics[width=0.47\textwidth,height=0.37\textwidth]{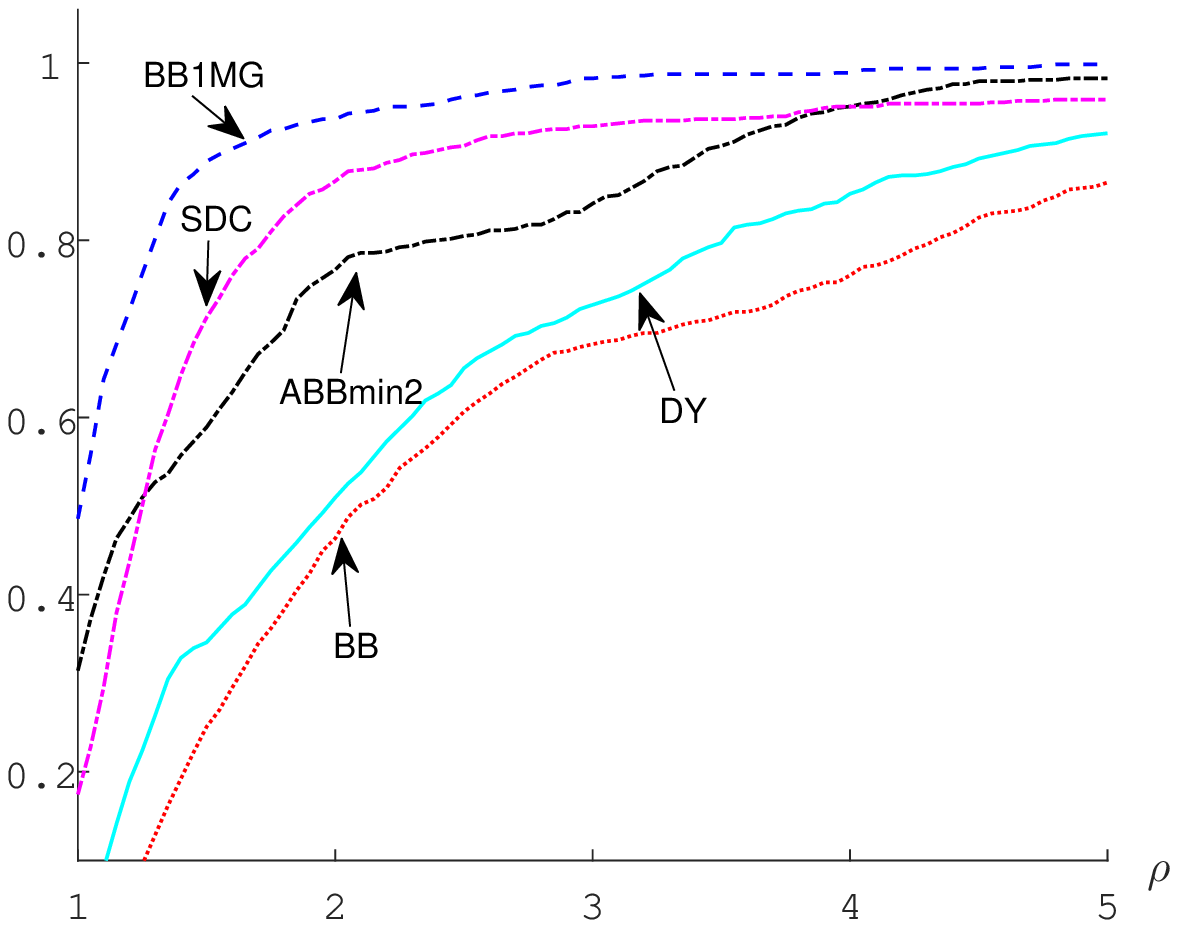}}
 \caption{Performance profiles for BB1SD (left)/BB1MG (right), and  BB, DY, ABBmin2 and SDC, iteration metric, 630 instances of the problems in Table \ref{tbspe}.}
\label{fig:fBB1}
\end{figure}

\begin{figure}[tbhp]
\centering
\subfloat{\label{fig:3a}\includegraphics[width=0.47\textwidth,height=0.37\textwidth]{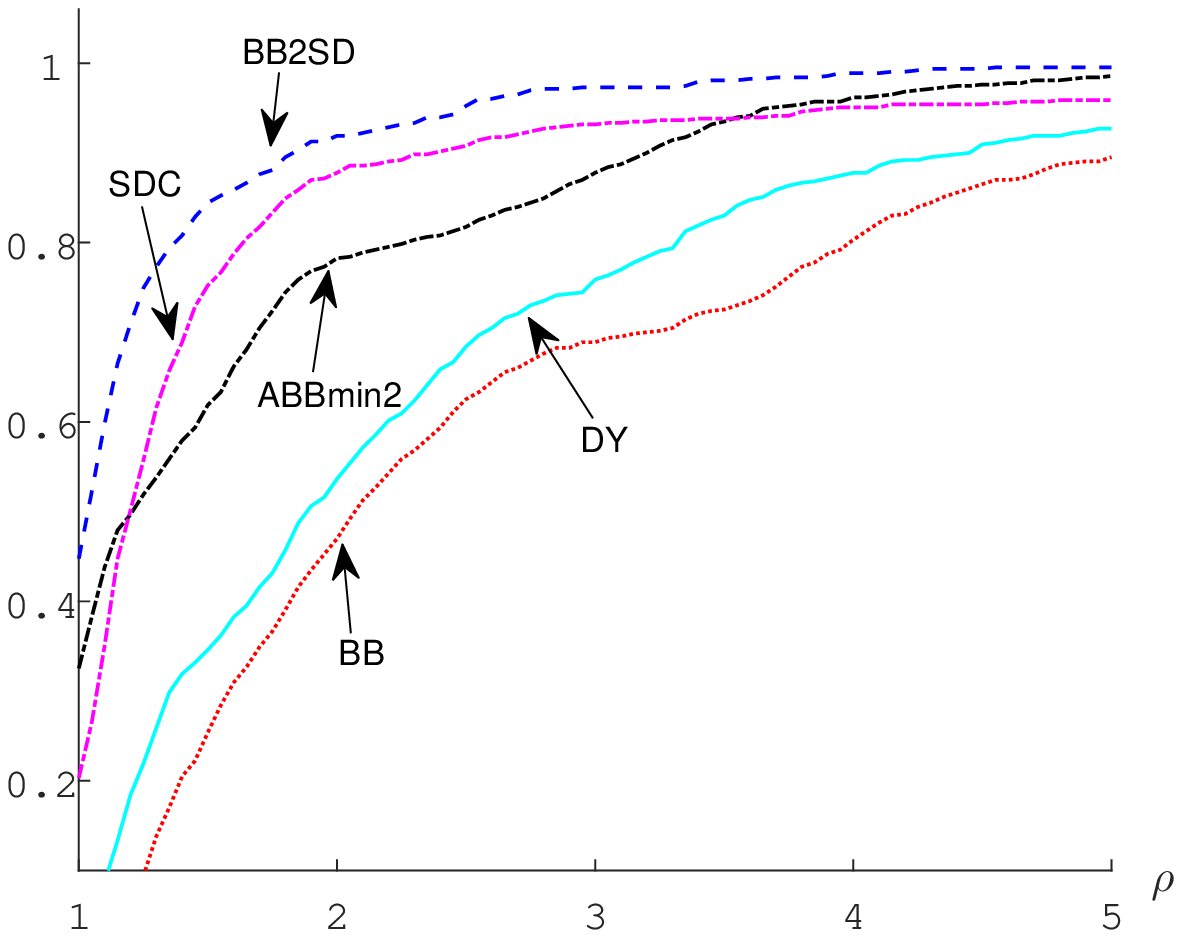}}
\subfloat{\label{fig:3b}\includegraphics[width=0.47\textwidth,height=0.37\textwidth]{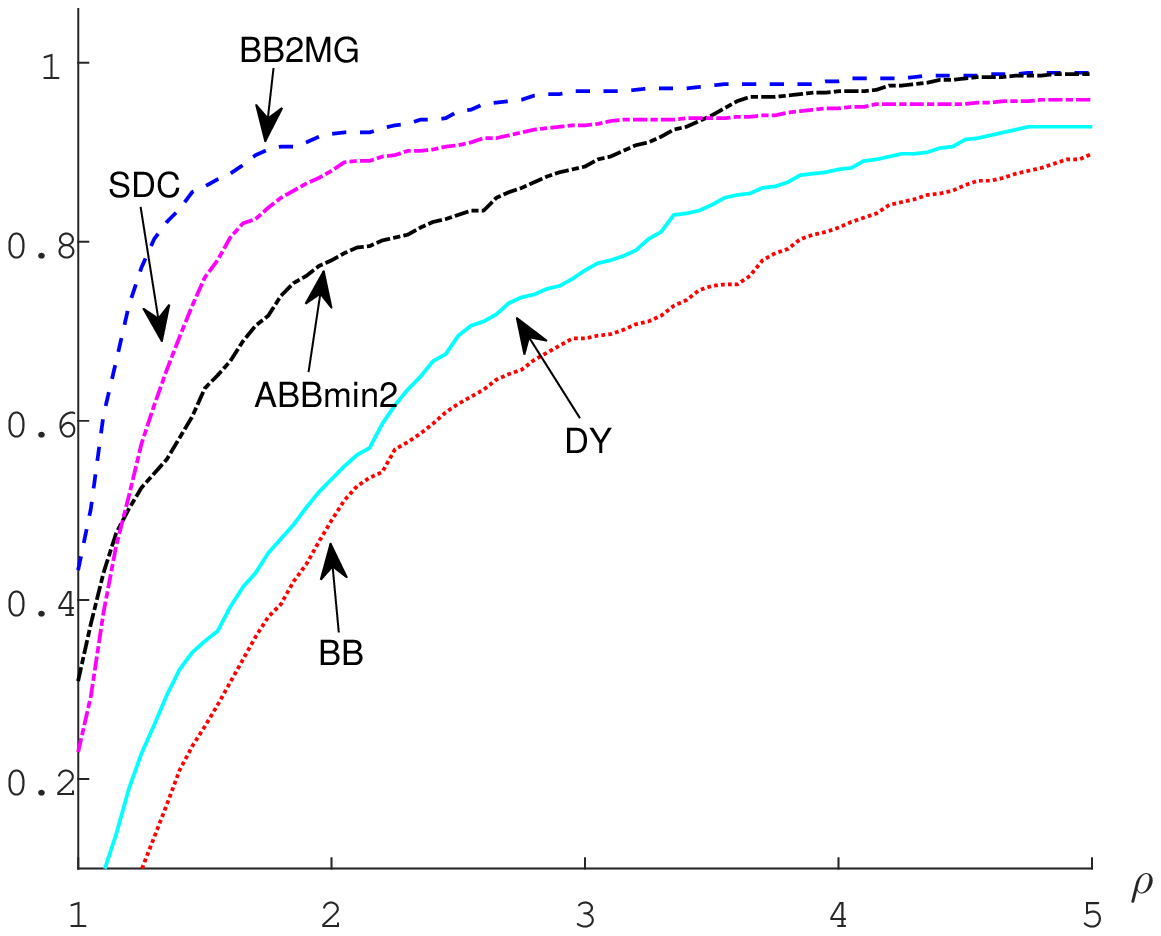}}
 \caption{Performance profiles for BB2SD (left)/BB2MG (right), and BB, DY, ABBmin2 and SDC, iteration metric, 630 instances of the problems in Table \ref{tbspe}.}
\label{fig:fBB2}
\end{figure}

\begin{figure}[ht!b]
\centering
    \includegraphics[width=0.5\textwidth,height=0.37\textwidth]{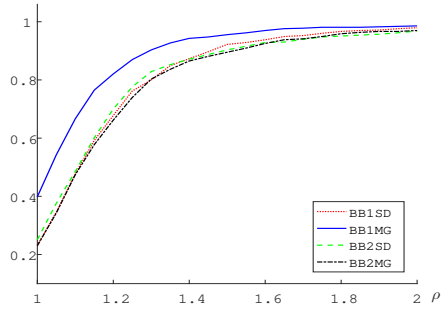}
  \caption{Performance profiles for BB1SD, BB1MG, BB2SD and BB2MG, iteration metric, 630 instances of the problems in Table \ref{tbspe}.}\label{fBBSDMG}
\end{figure}

We further compared these methods in Figures \ref{fig:fBB1} and \ref{fig:fBB2} by using the performance profiles of Dolan and Mor\'{e} \cite{dolan2002} on the iteration metric. In these figures, the vertical axis shows the percentage of the problems the method solves within the factor $\rho$ of the metric used by the most effective method in this comparison. We select the results of our four methods corresponding to the column $(15,15)$ in the above tables.
 It can be seen  that all our methods BB1SD, BB1MG, BB2SD and BB2MG clearly outperform the other compared methods.
For comparison of BB1SD, BB1MG, BB2SD and BB2MG, Figure \ref{fBBSDMG} shows that BB1MG is slightly better than the other three methods,
while  BB1SD, BB2SD and BB2MG do not show much difference in this test.

%\begin{figure}[ht!b]
%\centering
%    \includegraphics[width=0.52\textwidth,height=0.4\textwidth]{BB1SD-ar.pdf}
%  \caption{Performance profiles for BB1SD, BB, DY, ABBmin2 and SDC, iteration metric, 630 instances of the problems in Table \ref{tbspe}.}\label{fBB1SD}
%\end{figure}
%
%\begin{figure}[ht!b]
%\centering
%    \includegraphics[width=0.52\textwidth,height=0.4\textwidth]{BB1MG-ar.pdf}
%  \caption{Performance profiles for BB1MG, BB, DY, ABBmin2 and SDC, iteration metric, 630 instances of the problems in Table \ref{tbspe}.}\label{fBB1MG}
%\end{figure}

%\begin{figure}[ht!b]
%\centering
%    \includegraphics[width=0.52\textwidth,height=0.4\textwidth]{BB2SD-ar.pdf}
%\caption{Performance profiles for BB2SD, BB, DY, ABBmin2 and SDC, iteration metric, 630 instances of the problems in Table \ref{tbspe}.}\label{fBB2SD}
%\end{figure}
%
%\begin{figure}[ht!b]
%\centering
%    \includegraphics[width=0.52\textwidth,height=0.4\textwidth]{BB2MG-ar.pdf}
%   \caption{Performance profiles for BB2MG, BB, DY, ABBmin2 and SDC, iteration metric, 630 instances of the problems in Table \ref{tbspe}.}\label{fBB2MG}
%\end{figure}

\section{Conclusions and discussions}\label{seccls}
We present theoretical analyses on the asymptotic behaviors of a family of gradient methods whose stepsize is given by \eqref{glstep1}, which includes the steepest descent and minimal gradient methods as special cases. It is shown that each method in this family will asymptotically zigzag in a two-dimensional subspace spanned by the two eigenvectors corresponding to the largest and smallest eigenvalues of the Hessian. In order to accelerate the gradient methods,
we exploit the spectral property of a new stepsize to break the zigzagging pattern. This new stepsize is derived by imposing finite termination on minimizing two-dimensional strongly convex quadratics and is proved to converge to the reciprocal of the largest eigenvalue of
the Hessian for general $n$-dimensional case. Finally, we propose a very efficient periodic gradient method that alternately uses the BB stepsize,
 $\alpha_k$ in \eqref{glstep1} and our new stepsize. Our numerical results indicate that, by exploiting the asymptotic behavior and spectral properties of stepsizes, gradient methods can be greatly accelerated to outperform the BB method and other recently developed state-of-the-art gradient methods.
%We believe the theoretical results will be helpful in understanding and developing fast gradient methods.

As a final remark, one may also break the zigzagging pattern by employing the spectral property in \eqref{sum2reps}. In particular, we could use the following stepsize
  \begin{equation}\label{sum2reps2}
    \hat{\alpha}_k=\left(\frac{1}{\alpha_{2k}}+\frac{1}{\alpha_{2k+1}}\right)^{-1},
  \end{equation}
to break the zigzagging pattern.  By \eqref{sum2reps},  $ \hat{\alpha}_k$  satisfies
  \[
        \lim_{k\rightarrow\infty}\hat{\alpha}_k
    = \frac{1}{\lambda_{1}+\lambda_{n}}.
  \]
Hence, $ \hat{\alpha}_k$ is also a good approximation of $1/\lambda_n$ when the condition number $\kappa=\lambda_n/\lambda_1$ is large.
One may see the strategy used in \cite{de2013spectral} for the case of the SD method.

\appendix
\section{Tables}
%\lipsum[71]

\begin{table}[tbhp]
{\tiny
\caption{Number of averaged iterations of BB1SD, BB, DY, ABBmin2 and SDC on the problems in Table \ref{tbspe}.}\label{tbrandpBB1SD}
\begin{center}
\setlength{\tabcolsep}{0.48ex}
\begin{tabular}{|c|c|c|c|c|c|c|c|c|c|c|c|c|c|c|}
\hline
 \multicolumn{1}{|c|}{\multirow{2}{*}{Set}} &\multicolumn{1}{c|}{\multirow{2}{*}{$\epsilon$}}
 &\multicolumn{9}{c|}{$(K_m,K_s)$} &\multirow{2}{*}{BB} &\multirow{2}{*}{DY} &\multirow{2}{*}{ABBmin2} &\multirow{2}{*}{SDC}\\
 \cline{3-11}
 \multicolumn{1}{|c|}{}& \multicolumn{1}{c|}{}   &$(9,9)$ &$(9,13)$ &$(9,15)$ &$(13,9)$  &$(13,13)$  &$(13,15)$ &$(15,9)$ &$(15,13)$ &$(15,15)$
 & & & &\\
\hline

\multicolumn{1}{|c|}{\multirow{3}{*}{1}}
&\multicolumn{1}{c|}{\multirow{1}{*}{$10^{-6}$}}	   & 367.6   & 339.7   & 372.3 & 346.1   & 352.0   & 344.9 & 336.8   & 317.6   & 368.0    &\multirow{1}{*}{ 458.7} &\multirow{1}{*}{ 350.0} &\multirow{1}{*}{ 258.5} &\multirow{1}{*}{ 394.1}\\
&\multicolumn{1}{c|}{\multirow{1}{*}{$10^{-9}$}}	   &1232.3   & 983.7   &1312.7 &1149.2   &1149.5   &1281.4 &1011.4   &1086.7   &1150.6    &\multirow{1}{*}{3694.4} &\multirow{1}{*}{3520.9} &\multirow{1}{*}{ 511.2} &\multirow{1}{*}{2410.6}\\
&\multicolumn{1}{c|}{\multirow{1}{*}{$10^{-12}$}}    &1849.9   &1514.1   &1812.8 &1760.6   &1780.3   &1792.6 &1518.1   &1605.0   &1465.4    &\multirow{1}{*}{6825.4} &\multirow{1}{*}{6561.6} &\multirow{1}{*}{ 678.2} &\multirow{1}{*}{4917.2}\\
\hline

\multicolumn{1}{|c|}{\multirow{3}{*}{2}}
&\multicolumn{1}{c|}{\multirow{1}{*}{$10^{-6}$}}    & 242.1   & 244.4   & 235.8 & 238.2   & 229.3   & 236.7 & 249.6   & 233.9   & 240.9  &\multirow{1}{*}{ 455.7} &\multirow{1}{*}{ 406.7} &\multirow{1}{*}{ 380.0} &\multirow{1}{*}{ 234.1}\\
&\multicolumn{1}{c|}{\multirow{1}{*}{$10^{-9}$}}    & 816.3   & 790.9   & 765.5 & 840.0   & 729.8   & 737.0 & 758.9   & 750.3   & 746.8  &\multirow{1}{*}{1882.0} &\multirow{1}{*}{1682.6} &\multirow{1}{*}{1425.7} &\multirow{1}{*}{ 879.8}\\
&\multicolumn{1}{c|}{\multirow{1}{*}{$10^{-12}$}}   &1255.9   &1222.6   &1207.4 &1305.8   &1179.1   &1154.8 &1211.5   &1187.4   &1178.1  &\multirow{1}{*}{3149.5} &\multirow{1}{*}{2761.6} &\multirow{1}{*}{2255.9} &\multirow{1}{*}{1436.3}\\
\hline

\multicolumn{1}{|c|}{\multirow{3}{*}{3}}
&\multicolumn{1}{c|}{\multirow{1}{*}{$10^{-6}$}}    & 297.1   & 288.5   & 275.9 & 284.6   & 283.1   & 273.3 & 283.6   & 279.7   & 270.0    &\multirow{1}{*}{ 495.6} &\multirow{1}{*}{ 435.8} &\multirow{1}{*}{ 487.4} &\multirow{1}{*}{ 298.4}\\
&\multicolumn{1}{c|}{\multirow{1}{*}{$10^{-9}$}}    & 829.0   & 816.2   & 796.5 & 848.0   & 758.5   & 763.2 & 796.9   & 791.7   & 743.7    &\multirow{1}{*}{1859.9} &\multirow{1}{*}{1678.7} &\multirow{1}{*}{1509.2} &\multirow{1}{*}{ 926.1}\\
&\multicolumn{1}{c|}{\multirow{1}{*}{$10^{-12}$}}   &1330.5   &1241.0   &1252.8 &1345.9   &1224.3   &1176.1 &1275.2   &1189.9   &1178.9    &\multirow{1}{*}{3230.2} &\multirow{1}{*}{2747.0} &\multirow{1}{*}{2492.1} &\multirow{1}{*}{1402.4}\\
\hline

\multicolumn{1}{|c|}{\multirow{3}{*}{4}}
&\multicolumn{1}{c|}{\multirow{1}{*}{$10^{-6}$}}    & 358.0   & 331.8   & 343.5 & 331.6   & 331.6   & 318.4 & 331.5   & 326.6   & 342.6   &\multirow{1}{*}{ 715.0} &\multirow{1}{*}{ 585.0} &\multirow{1}{*}{ 679.9} &\multirow{1}{*}{ 345.7}\\
&\multicolumn{1}{c|}{\multirow{1}{*}{$10^{-9}$}}    & 882.6   & 823.2   & 825.3 & 917.8   & 814.2   & 817.4 & 860.3   & 808.6   & 832.4   &\multirow{1}{*}{2097.1} &\multirow{1}{*}{1927.2} &\multirow{1}{*}{1749.7} &\multirow{1}{*}{ 969.2}\\
&\multicolumn{1}{c|}{\multirow{1}{*}{$10^{-12}$}}   &1422.0   &1327.9   &1347.9 &1324.3   &1232.8   &1271.9 &1318.4   &1288.9   &1258.0   &\multirow{1}{*}{3355.7} &\multirow{1}{*}{3140.5} &\multirow{1}{*}{2673.9} &\multirow{1}{*}{1451.2}\\
\hline

\multicolumn{1}{|c|}{\multirow{3}{*}{5}}
&\multicolumn{1}{c|}{\multirow{1}{*}{$10^{-6}$}}    & 838.4   & 829.4   & 850.8 & 851.9   & 836.4   & 855.5 & 874.1   & 856.5   & 844.5        &\multirow{1}{*}{1091.5} &\multirow{1}{*}{ 849.1} &\multirow{1}{*}{1043.1} &\multirow{1}{*}{ 861.7}\\
&\multicolumn{1}{c|}{\multirow{1}{*}{$10^{-9}$}}    &3147.9   &3086.3   &2985.3 &2932.6   &3004.0   &3062.6 &3093.1   &3086.7   &3094.2        &\multirow{1}{*}{5262.6} &\multirow{1}{*}{4606.2} &\multirow{1}{*}{3542.8} &\multirow{1}{*}{4075.9}\\
&\multicolumn{1}{c|}{\multirow{1}{*}{$10^{-12}$}}   &4942.5   &4996.4   &4688.7 &4542.9   &5020.5   &4921.7 &4900.0   &4845.1   &4868.1        &\multirow{1}{*}{7803.1} &\multirow{1}{*}{8048.4} &\multirow{1}{*}{5518.2} &\multirow{1}{*}{6279.4}\\
\hline

\multicolumn{1}{|c|}{\multirow{3}{*}{6}}
&\multicolumn{1}{c|}{\multirow{1}{*}{$10^{-6}$}}    & 155.1   & 140.8   & 140.3 & 138.9   & 139.4   & 137.8 & 132.8   & 137.9   & 137.3     &\multirow{1}{*}{ 257.0} &\multirow{1}{*}{ 186.1} &\multirow{1}{*}{ 151.8} &\multirow{1}{*}{ 143.8}\\
&\multicolumn{1}{c|}{\multirow{1}{*}{$10^{-9}$}}    & 554.3   & 557.4   & 541.4 & 590.8   & 513.8   & 500.1 & 559.9   & 539.4   & 512.9     &\multirow{1}{*}{1574.2} &\multirow{1}{*}{1265.4} &\multirow{1}{*}{ 617.8} &\multirow{1}{*}{ 639.2}\\
&\multicolumn{1}{c|}{\multirow{1}{*}{$10^{-12}$}}   & 905.9   & 883.1   & 897.7 & 939.9   & 801.1   & 824.3 & 925.9   & 895.9   & 814.9     &\multirow{1}{*}{2603.9} &\multirow{1}{*}{2419.3} &\multirow{1}{*}{ 894.6} &\multirow{1}{*}{1129.3}\\
\hline

\multicolumn{1}{|c|}{\multirow{3}{*}{7}}
&\multicolumn{1}{c|}{\multirow{1}{*}{$10^{-6}$}}    & 455.6   & 437.0   & 430.8 & 457.9   & 432.9   & 424.0 & 445.8   & 411.3   & 424.8     &\multirow{1}{*}{ 893.7} &\multirow{1}{*}{ 800.3} &\multirow{1}{*}{ 772.7} &\multirow{1}{*}{ 470.5}\\
&\multicolumn{1}{c|}{\multirow{1}{*}{$10^{-9}$}}    & 905.8   & 876.0   & 828.2 & 922.4   & 870.5   & 869.8 & 925.6   & 851.0   & 859.6     &\multirow{1}{*}{2110.7} &\multirow{1}{*}{1868.1} &\multirow{1}{*}{1613.9} &\multirow{1}{*}{ 936.6}\\
&\multicolumn{1}{c|}{\multirow{1}{*}{$10^{-12}$}}   &1349.8   &1323.1   &1265.4 &1374.2   &1278.5   &1267.2 &1319.2   &1252.1   &1240.3     &\multirow{1}{*}{3252.1} &\multirow{1}{*}{2748.7} &\multirow{1}{*}{2372.9} &\multirow{1}{*}{1331.5}\\
\hline

\multicolumn{1}{|c|}{\multirow{3}{*}{total}}
&\multicolumn{1}{c|}{\multirow{1}{*}{$10^{-6}$}}  &2713.9 &2611.6 &2649.4 &2649.2 &2604.7 &2590.6 &2654.2 &{2563.5} &2628.1 &4367.2 &3613.0 &3773.4 &2748.3                 \\
&\multicolumn{1}{c|}{\multirow{1}{*}{$10^{-9}$}}  &8368.2 &7933.7 &8054.9 &8200.8 &{7840.3} &8031.5 &8006.1 &7914.4 &7940.2 &18480.9 &16549.1 &10970.3 &10837.4            \\
&\multicolumn{1}{c|}{\multirow{1}{*}{$10^{-12}$}} &13056.5 &12508.2 &12472.7 &12593.6 &12516.6 &12408.6 &12468.3 &12264.3 &{12003.7} &30219.9 &28427.1 &16885.8 &17947.3   \\
\hline

\end{tabular}
%\end{sideways}
\end{center}
}
\end{table}

\begin{table}[ht!b]
{\tiny
\caption{Number of averaged iterations of BB1MG, BB, DY, ABBmin2 and SDC on the problems in Table \ref{tbspe}.}
\begin{center}
\setlength{\tabcolsep}{0.48ex}
\begin{tabular}{|c|c|c|c|c|c|c|c|c|c|c|c|c|c|c|}
\hline
 \multicolumn{1}{|c|}{\multirow{2}{*}{Set}} &\multicolumn{1}{c|}{\multirow{2}{*}{$\epsilon$}}
 &\multicolumn{9}{c|}{$(K_m,K_s)$} &\multirow{2}{*}{BB} &\multirow{2}{*}{DY} &\multirow{2}{*}{ABBmin2} &\multirow{2}{*}{SDC}\\
 \cline{3-11}
 \multicolumn{1}{|c|}{}& \multicolumn{1}{c|}{}   &$(9,9)$ &$(9,13)$ &$(9,15)$ &$(13,9)$  &$(13,13)$  &$(13,15)$ &$(15,9)$ &$(15,13)$ &$(15,15)$
 & & & &\\
\hline

\multicolumn{1}{|c|}{\multirow{3}{*}{1}}
&\multicolumn{1}{c|}{\multirow{1}{*}{$10^{-6}$}}	 & 378.0   & 366.2   & 344.9 & 354.3   & 364.5   & 341.7 & 338.1   & 374.1   & 362.1       &\multirow{1}{*}{ 458.7} &\multirow{1}{*}{ 350.0} &\multirow{1}{*}{ 258.5} &\multirow{1}{*}{ 394.1}\\
&\multicolumn{1}{c|}{\multirow{1}{*}{$10^{-9}$}}	 &1187.6   &1369.2   &1192.8 &1029.0   &1297.6   &1040.6 &1124.6   &1201.2   &1095.8       &\multirow{1}{*}{3694.4} &\multirow{1}{*}{3520.9} &\multirow{1}{*}{ 511.2} &\multirow{1}{*}{2410.6}\\
&\multicolumn{1}{c|}{\multirow{1}{*}{$10^{-12}$}}  &1909.2   &1809.4   &1666.2 &1558.3   &1784.7   &1577.6 &1578.7   &1862.7   &1485.3       &\multirow{1}{*}{6825.4} &\multirow{1}{*}{6561.6} &\multirow{1}{*}{ 678.2} &\multirow{1}{*}{4917.2}\\
\hline

\multicolumn{1}{|c|}{\multirow{3}{*}{2}}
&\multicolumn{1}{c|}{\multirow{1}{*}{$10^{-6}$}}    & 216.5   & 211.0   & 227.0 & 218.2   & 211.2   & 228.5 & 223.3   & 225.5   & 230.2  &\multirow{1}{*}{ 455.7} &\multirow{1}{*}{ 406.7} &\multirow{1}{*}{ 380.0} &\multirow{1}{*}{ 234.1}\\
&\multicolumn{1}{c|}{\multirow{1}{*}{$10^{-9}$}}    & 729.7   & 679.9   & 703.0 & 665.9   & 674.4   & 686.0 & 675.6   & 665.8   & 680.9  &\multirow{1}{*}{1882.0} &\multirow{1}{*}{1682.6} &\multirow{1}{*}{1425.7} &\multirow{1}{*}{ 879.8}\\
&\multicolumn{1}{c|}{\multirow{1}{*}{$10^{-12}$}}   &1199.7   &1079.8   &1130.7 &1076.6   &1076.3   &1067.1 &1096.8   &1081.7   &1059.3  &\multirow{1}{*}{3149.5} &\multirow{1}{*}{2761.6} &\multirow{1}{*}{2255.9} &\multirow{1}{*}{1436.3}\\
\hline

\multicolumn{1}{|c|}{\multirow{3}{*}{3}}
&\multicolumn{1}{c|}{\multirow{1}{*}{$10^{-6}$}}    & 258.3   & 265.4   & 273.7 & 273.3   & 249.2   & 254.1 & 253.1   & 246.2   & 252.3   &\multirow{1}{*}{ 495.6} &\multirow{1}{*}{ 435.8} &\multirow{1}{*}{ 487.4} &\multirow{1}{*}{ 298.4}\\
&\multicolumn{1}{c|}{\multirow{1}{*}{$10^{-9}$}}    & 810.6   & 743.8   & 756.7 & 707.6   & 720.5   & 694.0 & 731.2   & 723.2   & 701.6   &\multirow{1}{*}{1859.9} &\multirow{1}{*}{1678.7} &\multirow{1}{*}{1509.2} &\multirow{1}{*}{ 926.1}\\
&\multicolumn{1}{c|}{\multirow{1}{*}{$10^{-12}$}}   &1208.6   &1137.4   &1182.6 &1112.4   &1128.5   &1102.7 &1153.6   &1108.9   &1099.7   &\multirow{1}{*}{3230.2} &\multirow{1}{*}{2747.0} &\multirow{1}{*}{2492.1} &\multirow{1}{*}{1402.4}\\
\hline

\multicolumn{1}{|c|}{\multirow{3}{*}{4}}
&\multicolumn{1}{c|}{\multirow{1}{*}{$10^{-6}$}}    & 309.8   & 325.1   & 305.4  & 315.2   & 309.3   & 312.6 & 315.1   & 304.9   & 315.8    &\multirow{1}{*}{ 715.0} &\multirow{1}{*}{ 585.0} &\multirow{1}{*}{ 679.9} &\multirow{1}{*}{ 345.7}\\
&\multicolumn{1}{c|}{\multirow{1}{*}{$10^{-9}$}}    & 871.1   & 753.9   & 764.6 & 771.7   & 766.2   & 748.4 & 766.8   & 749.2   & 768.3     &\multirow{1}{*}{2097.1} &\multirow{1}{*}{1927.2} &\multirow{1}{*}{1749.7} &\multirow{1}{*}{ 969.2}\\
&\multicolumn{1}{c|}{\multirow{1}{*}{$10^{-12}$}}   &1268.6   &1186.6   &1203.9 &1164.3   &1162.0   &1140.8 &1200.9   &1159.1   &1181.2     &\multirow{1}{*}{3355.7} &\multirow{1}{*}{3140.5} &\multirow{1}{*}{2673.9} &\multirow{1}{*}{1451.2}\\
\hline

\multicolumn{1}{|c|}{\multirow{3}{*}{5}}
&\multicolumn{1}{c|}{\multirow{1}{*}{$10^{-6}$}}    & 856.8   & 833.5   & 847.7 & 862.7   & 847.2   & 848.3 & 843.7   & 906.7   & 865.1       &\multirow{1}{*}{1091.5} &\multirow{1}{*}{ 849.1} &\multirow{1}{*}{1043.1} &\multirow{1}{*}{ 861.7}\\
&\multicolumn{1}{c|}{\multirow{1}{*}{$10^{-9}$}}    &3197.5   &3014.6   &3216.2 &2988.8   &3015.1   &3088.4 &3137.5   &3155.4   &3042.1       &\multirow{1}{*}{5262.6} &\multirow{1}{*}{4606.2} &\multirow{1}{*}{3542.8} &\multirow{1}{*}{4075.9}\\
&\multicolumn{1}{c|}{\multirow{1}{*}{$10^{-12}$}}   &4937.7   &4769.0   &4986.6 &4933.8   &4709.7   &4861.1 &4944.6   &5167.5   &4869.2       &\multirow{1}{*}{7803.1} &\multirow{1}{*}{8048.4} &\multirow{1}{*}{5518.2} &\multirow{1}{*}{6279.4}\\
\hline

\multicolumn{1}{|c|}{\multirow{3}{*}{6}}
&\multicolumn{1}{c|}{\multirow{1}{*}{$10^{-6}$}}    & 129.1   & 125.6   & 126.0 & 132.5   & 126.1   & 135.4 & 128.6   & 127.0   & 137.3     &\multirow{1}{*}{ 257.0} &\multirow{1}{*}{ 186.1} &\multirow{1}{*}{ 151.8} &\multirow{1}{*}{ 143.8}\\
&\multicolumn{1}{c|}{\multirow{1}{*}{$10^{-9}$}}    & 510.8   & 498.9   & 510.1 & 496.3   & 452.1   & 471.3 & 461.6   & 487.2   & 447.6     &\multirow{1}{*}{1574.2} &\multirow{1}{*}{1265.4} &\multirow{1}{*}{ 617.8} &\multirow{1}{*}{ 639.2}\\
&\multicolumn{1}{c|}{\multirow{1}{*}{$10^{-12}$}}   & 841.4   & 799.5   & 789.0 & 808.8   & 712.1   & 780.5 & 754.2   & 748.2   & 699.8     &\multirow{1}{*}{2603.9} &\multirow{1}{*}{2419.3} &\multirow{1}{*}{ 894.6} &\multirow{1}{*}{1129.3}\\
\hline

\multicolumn{1}{|c|}{\multirow{3}{*}{7}}
&\multicolumn{1}{c|}{\multirow{1}{*}{$10^{-6}$}}    & 400.6   & 417.1   & 382.8 & 423.1   & 407.0   & 405.6 & 402.0   & 415.8   & 402.7     &\multirow{1}{*}{ 893.7} &\multirow{1}{*}{ 800.3} &\multirow{1}{*}{ 772.7} &\multirow{1}{*}{ 470.5}\\
&\multicolumn{1}{c|}{\multirow{1}{*}{$10^{-9}$}}    & 841.3   & 815.6   & 788.3 & 832.9   & 820.8   & 794.4 & 825.4   & 844.7   & 814.5     &\multirow{1}{*}{2110.7} &\multirow{1}{*}{1868.1} &\multirow{1}{*}{1613.9} &\multirow{1}{*}{ 936.6}\\
&\multicolumn{1}{c|}{\multirow{1}{*}{$10^{-12}$}}   &1245.0   &1193.1   &1161.9 &1218.1   &1202.7   &1190.3 &1210.3   &1238.0   &1167.7     &\multirow{1}{*}{3252.1} &\multirow{1}{*}{2748.7} &\multirow{1}{*}{2372.9} &\multirow{1}{*}{1331.5}\\
\hline

\multicolumn{1}{|c|}{\multirow{3}{*}{total}}
&\multicolumn{1}{c|}{\multirow{1}{*}{$10^{-6}$}}    &2549.1 &2543.9 &2507.5 &2579.3 &2514.5 &2526.2 &{2503.9} &2600.2 &2565.5 &4367.2 &3613.0 &3773.4 &2748.3                \\
&\multicolumn{1}{c|}{\multirow{1}{*}{$10^{-9}$}}    &8148.6 &7875.9 &7931.7 &{7492.2} &7746.7 &7523.1 &7722.7 &7826.7 &7550.8 &18480.9 &16549.1 &10970.3 &10837.4           \\
&\multicolumn{1}{c|}{\multirow{1}{*}{$10^{-12}$}}   &12610.2 &11974.8 &12120.9 &11872.3 &11776.0 &11720.1 &11939.1 &12366.1 &{11562.2} &30219.9 &28427.1 &16885.8 &17947.3  \\
\hline

\end{tabular}
\end{center}
}\label{tbrandpBB1MG}
\end{table}

\begin{table}[ht!b]
{\tiny
\caption{Number of averaged iterations of BB2SD, BB, DY, ABBmin2 and SDC on the problems in Table \ref{tbspe}.}
\begin{center}
\setlength{\tabcolsep}{0.48ex}
\begin{tabular}{|c|c|c|c|c|c|c|c|c|c|c|c|c|c|c|}
\hline
 \multicolumn{1}{|c|}{\multirow{2}{*}{Set}} &\multicolumn{1}{c|}{\multirow{2}{*}{$\epsilon$}}
 &\multicolumn{9}{c|}{$(K_m,K_s)$} &\multirow{2}{*}{BB} &\multirow{2}{*}{DY} &\multirow{2}{*}{ABBmin2} &\multirow{2}{*}{SDC}\\
 \cline{3-11}
 \multicolumn{1}{|c|}{}& \multicolumn{1}{c|}{}   &$(9,9)$ &$(9,13)$ &$(9,15)$ &$(13,9)$  &$(13,13)$  &$(13,15)$ &$(15,9)$ &$(15,13)$ &$(15,15)$
 & & & &\\
\hline

\multicolumn{1}{|c|}{\multirow{3}{*}{1}}
&\multicolumn{1}{c|}{\multirow{1}{*}{$10^{-6}$}}	 & 347.9   & 357.2   & 365.1 & 349.4   & 344.3   & 325.0 & 338.1   & 349.4   & 369.2      &\multirow{1}{*}{ 458.7} &\multirow{1}{*}{ 350.0} &\multirow{1}{*}{ 258.5} &\multirow{1}{*}{ 394.1}\\
&\multicolumn{1}{c|}{\multirow{1}{*}{$10^{-9}$}}	 &1132.2   &1454.1   &1247.4 &1192.7   &1224.4   &1274.7 &1237.7   &1291.9   &1209.6      &\multirow{1}{*}{3694.4} &\multirow{1}{*}{3520.9} &\multirow{1}{*}{ 511.2} &\multirow{1}{*}{2410.6}\\
&\multicolumn{1}{c|}{\multirow{1}{*}{$10^{-12}$}}  &1985.3   &2429.8   &1986.8 &1838.2   &2062.1   &2181.2 &1958.2   &1961.0   &1927.2      &\multirow{1}{*}{6825.4} &\multirow{1}{*}{6561.6} &\multirow{1}{*}{ 678.2} &\multirow{1}{*}{4917.2}\\
\hline

\multicolumn{1}{|c|}{\multirow{3}{*}{2}}
&\multicolumn{1}{c|}{\multirow{1}{*}{$10^{-6}$}}    & 219.4   & 223.9   & 220.5 & 226.0   & 229.3   & 224.4 & 217.8   & 220.4   & 226.4  &\multirow{1}{*}{ 455.7} &\multirow{1}{*}{ 406.7} &\multirow{1}{*}{ 380.0} &\multirow{1}{*}{ 234.1}\\
&\multicolumn{1}{c|}{\multirow{1}{*}{$10^{-9}$}}    & 749.4   & 723.3   & 713.9 & 746.6   & 720.1   & 711.2 & 728.1   & 729.4   & 713.3  &\multirow{1}{*}{1882.0} &\multirow{1}{*}{1682.6} &\multirow{1}{*}{1425.7} &\multirow{1}{*}{ 879.8}\\
&\multicolumn{1}{c|}{\multirow{1}{*}{$10^{-12}$}}   &1235.9   &1188.4   &1168.4 &1167.9   &1158.1   &1158.3 &1165.2   &1186.0   &1130.9  &\multirow{1}{*}{3149.5} &\multirow{1}{*}{2761.6} &\multirow{1}{*}{2255.9} &\multirow{1}{*}{1436.3}\\
\hline

\multicolumn{1}{|c|}{\multirow{3}{*}{3}}
&\multicolumn{1}{c|}{\multirow{1}{*}{$10^{-6}$}}    & 248.5   & 259.0   & 253.8 & 254.0   & 246.3   & 261.6 & 252.6   & 262.8   & 267.4   &\multirow{1}{*}{ 495.6} &\multirow{1}{*}{ 435.8} &\multirow{1}{*}{ 487.4} &\multirow{1}{*}{ 298.4}\\
&\multicolumn{1}{c|}{\multirow{1}{*}{$10^{-9}$}}    & 780.5   & 757.1   & 754.2 & 759.3   & 738.4   & 767.2 & 793.6   & 774.4   & 759.3   &\multirow{1}{*}{1859.9} &\multirow{1}{*}{1678.7} &\multirow{1}{*}{1509.2} &\multirow{1}{*}{ 926.1}\\
&\multicolumn{1}{c|}{\multirow{1}{*}{$10^{-12}$}}   &1229.4   &1230.7   &1227.8 &1216.0   &1214.8   &1182.3 &1215.2   &1227.7   &1210.6   &\multirow{1}{*}{3230.2} &\multirow{1}{*}{2747.0} &\multirow{1}{*}{2492.1} &\multirow{1}{*}{1402.4}\\
\hline

\multicolumn{1}{|c|}{\multirow{3}{*}{4}}
&\multicolumn{1}{c|}{\multirow{1}{*}{$10^{-6}$}}    & 320.8   & 315.1   & 305.5 & 313.6   & 315.9   & 310.9 & 318.4   & 307.5   & 317.1    &\multirow{1}{*}{ 715.0} &\multirow{1}{*}{ 585.0} &\multirow{1}{*}{ 679.9} &\multirow{1}{*}{ 345.7}\\
&\multicolumn{1}{c|}{\multirow{1}{*}{$10^{-9}$}}    & 805.0   & 823.3   & 813.4 & 819.5   & 813.5   & 789.0 & 779.5   & 836.1   & 802.5    &\multirow{1}{*}{2097.1} &\multirow{1}{*}{1927.2} &\multirow{1}{*}{1749.7} &\multirow{1}{*}{ 969.2}\\
&\multicolumn{1}{c|}{\multirow{1}{*}{$10^{-12}$}}   &1348.7   &1298.3   &1244.4 &1242.8   &1276.1   &1238.6 &1250.0   &1269.9   &1246.3    &\multirow{1}{*}{3355.7} &\multirow{1}{*}{3140.5} &\multirow{1}{*}{2673.9} &\multirow{1}{*}{1451.2}\\
\hline

\multicolumn{1}{|c|}{\multirow{3}{*}{5}}
&\multicolumn{1}{c|}{\multirow{1}{*}{$10^{-6}$}}    & 860.0   & 847.3   & 848.7 & 831.2   & 799.3   & 825.5 & 804.4   & 809.5   & 862.0       &\multirow{1}{*}{1091.5} &\multirow{1}{*}{ 849.1} &\multirow{1}{*}{1043.1} &\multirow{1}{*}{ 861.7}\\
&\multicolumn{1}{c|}{\multirow{1}{*}{$10^{-9}$}}    &3066.6   &3191.0   &2998.8 &2918.1   &3049.0   &3038.7 &2995.5   &2995.7   &3095.7       &\multirow{1}{*}{5262.6} &\multirow{1}{*}{4606.2} &\multirow{1}{*}{3542.8} &\multirow{1}{*}{4075.9}\\
&\multicolumn{1}{c|}{\multirow{1}{*}{$10^{-12}$}}   &5272.4   &5133.8   &5106.8 &4962.9   &4867.3   &4894.1 &5083.6   &4775.5   &5100.4       &\multirow{1}{*}{7803.1} &\multirow{1}{*}{8048.4} &\multirow{1}{*}{5518.2} &\multirow{1}{*}{6279.4}\\
\hline

\multicolumn{1}{|c|}{\multirow{3}{*}{6}}
&\multicolumn{1}{c|}{\multirow{1}{*}{$10^{-6}$}}    & 129.1   & 138.8   & 124.8 & 128.4   & 135.3   & 133.7 & 122.2   & 130.8   & 133.4     &\multirow{1}{*}{ 257.0} &\multirow{1}{*}{ 186.1} &\multirow{1}{*}{ 151.8} &\multirow{1}{*}{ 143.8}\\
&\multicolumn{1}{c|}{\multirow{1}{*}{$10^{-9}$}}    & 560.3   & 549.5   & 531.5 & 514.6   & 520.9   & 538.9 & 516.5   & 530.4   & 525.1     &\multirow{1}{*}{1574.2} &\multirow{1}{*}{1265.4} &\multirow{1}{*}{ 617.8} &\multirow{1}{*}{ 639.2}\\
&\multicolumn{1}{c|}{\multirow{1}{*}{$10^{-12}$}}   & 912.8   & 892.1   & 940.0 & 913.5   & 928.1   & 873.5 & 892.5   & 873.3   & 845.2     &\multirow{1}{*}{2603.9} &\multirow{1}{*}{2419.3} &\multirow{1}{*}{ 894.6} &\multirow{1}{*}{1129.3}\\
\hline

\multicolumn{1}{|c|}{\multirow{3}{*}{7}}
&\multicolumn{1}{c|}{\multirow{1}{*}{$10^{-6}$}}    & 418.4   & 393.6   & 406.6 & 410.6   & 409.7   & 418.4 & 394.8   & 429.7   & 405.9     &\multirow{1}{*}{ 893.7} &\multirow{1}{*}{ 800.3} &\multirow{1}{*}{ 772.7} &\multirow{1}{*}{ 470.5}\\
&\multicolumn{1}{c|}{\multirow{1}{*}{$10^{-9}$}}    & 898.0   & 835.8   & 849.0 & 852.9   & 847.6   & 847.8 & 868.4   & 873.3   & 848.4     &\multirow{1}{*}{2110.7} &\multirow{1}{*}{1868.1} &\multirow{1}{*}{1613.9} &\multirow{1}{*}{ 936.6}\\
&\multicolumn{1}{c|}{\multirow{1}{*}{$10^{-12}$}}   &1324.7   &1238.8   &1221.1 &1290.1   &1263.3   &1265.1 &1302.7   &1279.2   &1267.4     &\multirow{1}{*}{3252.1} &\multirow{1}{*}{2748.7} &\multirow{1}{*}{2372.9} &\multirow{1}{*}{1331.5}\\
\hline

\multicolumn{1}{|c|}{\multirow{3}{*}{total}}
&\multicolumn{1}{c|}{\multirow{1}{*}{$10^{-6}$}}     &2544.1 &2534.9 &2525.0 &2513.2 &2480.1 &2499.5 &{2448.3} &2510.1 &2581.4 &4367.2 &3613.0 &3773.4 &2748.3                \\
&\multicolumn{1}{c|}{\multirow{1}{*}{$10^{-9}$}}     &7992.0 &8334.1 &7908.2 &{7803.7} &7913.9 &7967.5 &7919.3 &8031.2 &7953.9 &18480.9 &16549.1 &10970.3 &10837.4           \\
&\multicolumn{1}{c|}{\multirow{1}{*}{$10^{-12}$}}    &13309.2 &13411.9 &12895.3 &12631.4 &12769.8 &12793.1 &12867.4 &{12572.6} &12728.0 &30219.9 &28427.1 &16885.8 &17947.3  \\
\hline

\end{tabular}
\end{center}
}\label{tbrandpBB2SD}
\end{table}

\begin{table}[ht!b]
{\tiny
\caption{Number of averaged iterations of BB2MG, BB, DY, ABBmin2 and SDC on the problems in Table \ref{tbspe}.}
\begin{center}
\setlength{\tabcolsep}{0.48ex}
\begin{tabular}{|c|c|c|c|c|c|c|c|c|c|c|c|c|c|c|}
\hline
 \multicolumn{1}{|c|}{\multirow{2}{*}{Set}} &\multicolumn{1}{c|}{\multirow{2}{*}{$\epsilon$}}
 &\multicolumn{9}{c|}{$(K_m,K_s)$} &\multirow{2}{*}{BB} &\multirow{2}{*}{DY} &\multirow{2}{*}{ABBmin2} &\multirow{2}{*}{SDC}\\
 \cline{3-11}
 \multicolumn{1}{|c|}{}& \multicolumn{1}{c|}{}   &$(9,9)$ &$(9,13)$ &$(9,15)$ &$(13,9)$  &$(13,13)$  &$(13,15)$ &$(15,9)$ &$(15,13)$ &$(15,15)$
 & & & &\\
\hline

\multicolumn{1}{|c|}{\multirow{3}{*}{1}}
&\multicolumn{1}{c|}{\multirow{1}{*}{$10^{-6}$}}	 & 355.7   & 365.1   & 341.9 & 322.6   & 350.7   & 327.5  & 337.9   & 313.6   & 321.1       &\multirow{1}{*}{ 458.7} &\multirow{1}{*}{ 350.0} &\multirow{1}{*}{ 258.5} &\multirow{1}{*}{ 394.1}\\
&\multicolumn{1}{c|}{\multirow{1}{*}{$10^{-9}$}}	 &1209.5   &1327.4   & 908.0 &1064.7   &1206.9   &1209.7 & 965.6   &1255.1   &1351.1        &\multirow{1}{*}{3694.4} &\multirow{1}{*}{3520.9} &\multirow{1}{*}{ 511.2} &\multirow{1}{*}{2410.6}\\
&\multicolumn{1}{c|}{\multirow{1}{*}{$10^{-12}$}}  &1858.7   &1772.7   &1477.3  &1640.8   &1701.6   &1877.9  &1651.6   &1889.2   &1751.7      &\multirow{1}{*}{6825.4} &\multirow{1}{*}{6561.6} &\multirow{1}{*}{ 678.2} &\multirow{1}{*}{4917.2}\\
\hline

\multicolumn{1}{|c|}{\multirow{3}{*}{2}}
&\multicolumn{1}{c|}{\multirow{1}{*}{$10^{-6}$}}    & 235.1   & 237.9   & 238.2  & 233.0   & 229.2   & 239.2 & 236.4   & 235.2   & 238.0   &\multirow{1}{*}{ 455.7} &\multirow{1}{*}{ 406.7} &\multirow{1}{*}{ 380.0} &\multirow{1}{*}{ 234.1}\\
&\multicolumn{1}{c|}{\multirow{1}{*}{$10^{-9}$}}    & 822.7   & 778.9   & 752.8  & 805.0   & 747.0   & 762.7  & 785.7   & 748.0   & 737.0  &\multirow{1}{*}{1882.0} &\multirow{1}{*}{1682.6} &\multirow{1}{*}{1425.7} &\multirow{1}{*}{ 879.8}\\
&\multicolumn{1}{c|}{\multirow{1}{*}{$10^{-12}$}}   &1273.8   &1233.0   &1212.6 &1294.3   &1144.2   &1193.2 &1248.0   &1178.3   &1167.0    &\multirow{1}{*}{3149.5} &\multirow{1}{*}{2761.6} &\multirow{1}{*}{2255.9} &\multirow{1}{*}{1436.3}\\
\hline

\multicolumn{1}{|c|}{\multirow{3}{*}{3}}
&\multicolumn{1}{c|}{\multirow{1}{*}{$10^{-6}$}}    & 273.8   & 265.6   & 287.9 & 264.6   & 271.2   & 274.4 & 275.2   & 263.1   & 281.9     &\multirow{1}{*}{ 495.6} &\multirow{1}{*}{ 435.8} &\multirow{1}{*}{ 487.4} &\multirow{1}{*}{ 298.4}\\
&\multicolumn{1}{c|}{\multirow{1}{*}{$10^{-9}$}}    & 866.7   & 831.4   & 793.5  & 862.2   & 777.6   & 789.1  & 804.3   & 786.0   & 786.6   &\multirow{1}{*}{1859.9} &\multirow{1}{*}{1678.7} &\multirow{1}{*}{1509.2} &\multirow{1}{*}{ 926.1}\\
&\multicolumn{1}{c|}{\multirow{1}{*}{$10^{-12}$}}   &1313.6   &1318.9   &1244.3 &1361.4   &1219.6   &1234.7 &1313.4   &1271.2   &1251.8     &\multirow{1}{*}{3230.2} &\multirow{1}{*}{2747.0} &\multirow{1}{*}{2492.1} &\multirow{1}{*}{1402.4}\\
\hline

\multicolumn{1}{|c|}{\multirow{3}{*}{4}}
&\multicolumn{1}{c|}{\multirow{1}{*}{$10^{-6}$}}    & 333.7   & 335.8   & 341.9  & 353.0   & 319.9   & 317.4 & 331.7   & 333.0   & 329.1    &\multirow{1}{*}{ 715.0} &\multirow{1}{*}{ 585.0} &\multirow{1}{*}{ 679.9} &\multirow{1}{*}{ 345.7}\\
&\multicolumn{1}{c|}{\multirow{1}{*}{$10^{-9}$}}    & 876.9   & 877.7   & 853.3  & 863.8   & 844.5   & 836.6 & 881.4   & 804.5   & 800.1    &\multirow{1}{*}{2097.1} &\multirow{1}{*}{1927.2} &\multirow{1}{*}{1749.7} &\multirow{1}{*}{ 969.2}\\
&\multicolumn{1}{c|}{\multirow{1}{*}{$10^{-12}$}}   &1364.3   &1329.9   &1307.0  &1351.1   &1296.9   &1259.4 &1337.0   &1275.7   &1286.7    &\multirow{1}{*}{3355.7} &\multirow{1}{*}{3140.5} &\multirow{1}{*}{2673.9} &\multirow{1}{*}{1451.2}\\
\hline

\multicolumn{1}{|c|}{\multirow{3}{*}{5}}
&\multicolumn{1}{c|}{\multirow{1}{*}{$10^{-6}$}}    & 806.4   & 836.7   & 837.7 & 807.1   & 842.2   & 862.9 & 817.8   & 814.9   & 819.9     &\multirow{1}{*}{1091.5} &\multirow{1}{*}{ 849.1} &\multirow{1}{*}{1043.1} &\multirow{1}{*}{ 861.7}\\
&\multicolumn{1}{c|}{\multirow{1}{*}{$10^{-9}$}}    &3106.8   &3101.1   &3008.3 &3102.0   &3169.6   &3058.9 &3073.8   &2997.6   &3097.9     &\multirow{1}{*}{5262.6} &\multirow{1}{*}{4606.2} &\multirow{1}{*}{3542.8} &\multirow{1}{*}{4075.9}\\
&\multicolumn{1}{c|}{\multirow{1}{*}{$10^{-12}$}}   &4996.6   &5100.9   &4749.5 &5079.1   &5012.9   &5004.8 &5090.7   &5094.0   &4708.6     &\multirow{1}{*}{7803.1} &\multirow{1}{*}{8048.4} &\multirow{1}{*}{5518.2} &\multirow{1}{*}{6279.4}\\
\hline

\multicolumn{1}{|c|}{\multirow{3}{*}{6}}
&\multicolumn{1}{c|}{\multirow{1}{*}{$10^{-6}$}}    & 137.1   & 138.9   & 135.9 & 143.4   & 135.1   & 139.0 & 135.1   & 136.9   & 138.9       &\multirow{1}{*}{ 257.0} &\multirow{1}{*}{ 186.1} &\multirow{1}{*}{ 151.8} &\multirow{1}{*}{ 143.8}\\
&\multicolumn{1}{c|}{\multirow{1}{*}{$10^{-9}$}}    & 612.6   & 571.2   & 535.3 & 588.6   & 543.6   & 523.8 & 504.2   & 569.0   & 523.0       &\multirow{1}{*}{1574.2} &\multirow{1}{*}{1265.4} &\multirow{1}{*}{ 617.8} &\multirow{1}{*}{ 639.2}\\
&\multicolumn{1}{c|}{\multirow{1}{*}{$10^{-12}$}}   & 933.9   & 874.6   & 870.0  &1026.1   & 864.7   & 830.9  & 862.3   & 910.9   & 861.2     &\multirow{1}{*}{2603.9} &\multirow{1}{*}{2419.3} &\multirow{1}{*}{ 894.6} &\multirow{1}{*}{1129.3}\\
\hline

\multicolumn{1}{|c|}{\multirow{3}{*}{7}}
&\multicolumn{1}{c|}{\multirow{1}{*}{$10^{-6}$}}    & 462.7   & 430.8   & 434.4 & 454.2   & 428.2   & 438.8 & 440.8   & 437.9   & 435.1     &\multirow{1}{*}{ 893.7} &\multirow{1}{*}{ 800.3} &\multirow{1}{*}{ 772.7} &\multirow{1}{*}{ 470.5}\\
&\multicolumn{1}{c|}{\multirow{1}{*}{$10^{-9}$}}    & 957.1   & 932.7   & 904.4 & 935.3   & 868.1   & 889.4 & 933.5   & 917.1   & 869.6     &\multirow{1}{*}{2110.7} &\multirow{1}{*}{1868.1} &\multirow{1}{*}{1613.9} &\multirow{1}{*}{ 936.6}\\
&\multicolumn{1}{c|}{\multirow{1}{*}{$10^{-12}$}}   &1383.7   &1337.3   &1281.5 &1344.8   &1288.7   &1323.3 &1373.0   &1310.1   &1277.8     &\multirow{1}{*}{3252.1} &\multirow{1}{*}{2748.7} &\multirow{1}{*}{2372.9} &\multirow{1}{*}{1331.5}\\
\hline

\multicolumn{1}{|c|}{\multirow{3}{*}{total}}
&\multicolumn{1}{c|}{\multirow{1}{*}{$10^{-6}$}}     &2604.5 &2610.8 &2617.9 &2577.9 &2576.5 &2599.2 &2574.9 &{2534.6} &2564.0 &4367.2 &3613.0 &3773.4 &2748.3               \\
&\multicolumn{1}{c|}{\multirow{1}{*}{$10^{-9}$}}     &8452.3 &8420.4 &{7755.6} &8221.6 &8157.3 &8070.2 &7948.5 &8077.3 &8165.3 &18480.9 &16549.1 &10970.3 &10837.4          \\
&\multicolumn{1}{c|}{\multirow{1}{*}{$10^{-12}$}}    &13124.6 &12967.3 &{12142.2} &13097.6 &12528.6 &12724.2 &12876.0 &12929.4 &12304.8 &30219.9 &28427.1 &16885.8 &17947.3 \\
\hline
\end{tabular}
\end{center}
}\label{tbrandpBB2MG}
\end{table}

%\section*{Acknowledgments}
%We would like to acknowledge the assistance of volunteers in putting
%together this example manuscript and supplement.

%\bibliographystyle{siamplain}
%\bibliography{BB_hyk}

\end{document}